
\input amstex
\documentstyle{amsppt}
\NoBlackBoxes
\mag1100
\let\aro=@ 
\newcount\numcount
\def\numerote{\global\advance\numcount by 1 \the\numcount}
\def\lastnum[#1]{{\advance \numcount by #1(\the\numcount)}}
\topmatter
\title 
Foliated Structure of The Kuranishi Space and Isomorphisms of Deformation Families of Compact Complex Manifolds 
(Structure feuillet\'ee de l'espace de Kuranishi et isomorphismes de familles de d\'eformations de vari\'et\'es compactes complexes)
\endtitle 
\rightheadtext{Isomorphisms of Deformation Families}
\author
Laurent Meersseman
\endauthor
\date October 8th, 2010\enddate 
\address 
{Laurent Meersseman}\hfill\hfill\linebreak
\indent{I.M.B.}\hfill\hfill\linebreak
\indent{Universit\'e de Bourgogne}\hfill\hfill\linebreak
\indent{B.P. 47870}\hfill\hfill\linebreak
\indent{21078 Dijon Cedex}\hfill\hfill\linebreak
\indent{France}\hfill\hfill
\endaddress
\email laurent.meersseman\@u-bourgogne.fr \endemail

  
\thanks
This work was elaborated and written during a one-year visit as a CNRS member at the Pacific Institute for Mathematical Sciences and
the University of British Columbia, Vancouver, BC. I would like to thank the PIMS and the UBC Math Department for their hospitality, as
well as the CNRS for giving me the opportunity of this visit.
\medskip 
This work was partially supported by grant FABER from the Conseil R\'egional de Bourgogne and by project COMPLEXE (ANR-08-JCJC-0130-01) from the Agence Nationale de la Recherche. 
\endthanks

\abstract
Consider the following uniformization problem. Take two holomorphic
(parametrized by some analytic set defined on a neighborhood of $0$ in $\Bbb C^p$, for some $p>0$) or differentiable (parametrized by 
an open neighborhood of $0$ in $\Bbb R^p$, for some $p>0$) deformation families of compact complex manifolds. Assume they are pointwise
isomorphic, that is for each point $t$ of the parameter space, the fiber over $t$ of the
first family is biholomorphic to the fiber over $t$ of the second family. Then, under which
conditions are the two families locally isomorphic at 0?
In this article, we give a
sufficient condition in the case of holomorphic families. We show then that,
surprisingly, this condition is not sufficient in the case of differentiable families.
We also describe different types of counterexamples and give some elements of classification of the counterexamples.
These results rely on a geometric study of the Kuranishi space of a compact complex
manifold.

\medskip
Consid\'erons le probl\`eme d'uniformisation suivant. Prenons deux familles de d\'eformation holomorphes (param\'etr\'ees par un ensemble analytique d\'efini dans un voisinage de $0$ dans $\Bbb C^p$ pour $p>0$) ou
diff\'erentiables (param\'etr\'ees par un voisinage de $0$ dans $\Bbb R^p$ pour $p>0$) de vari\'et\'es compactes complexes. Supposons les ponctuellement isomorphes, c'est-\`a-dire, que pour tout point $t$
de l'espace des param\`etres, la fibre en $t$ de la premi\`ere famille est biholomorphe \`a la fibre en $t$ de la deuxi\`eme famille. Sous quelle(s) condition(s) les deux familles sont-elles localement isomorphes
en $0$? Dans cet article, nous donnons une condition suffisante dans le cas de familles holomorphes. Nous montrons ensuite que, de fa\c con surprenante, la condition n'est pas suffisante dans le cas des familles
diff\'erentiables. Nous d\'ecrivons \'egalement plusieurs types de contre-exemples et donnons quelques \'el\'ements de classifications de ces contre-exemples. Ces r\'esultats reposent sur une \'etude g\'eom\'etrique
de l'espace de Kuranishi d'une vari\'et\'e compacte complexe.
\endabstract

\endtopmatter

\document
\head
{\bf Introduction}
\endhead

This article deals with the problem of giving a useful criterion to ensure that two holomorphic (respectively differentiable) 
deformation families are isomorphic as families. This takes the form of the following uniformization problem. Let 
$$
\pi_i\ :\ \Cal X_i\to U \quad\text{ respectively }\quad \pi_i\ :\ \Cal X_i\to V
\leqno i=1,2
$$
be two holomorphic (respectively differentiable) families of compact complex manifolds parametrized
by some analytic set $U$ defined on a neighborhood of $0$ in $\Bbb C^p$, for some $p>0$ (respectively 
an open neighborhood $V$ of $0$ in $\Bbb R^p$, for some $p>0$). Assume that they are {\it pointwise isomorphic},
that is, for all $t\in U$ (respectively $t\in V$), the fiber $X_1(t)=\pi_1^{-1}(\{t\})$ is biholomorphic to the fiber
$X_2(t)=\pi_2^{-1}(\{t\})$. Then the question is

\proclaim{Question 1} Under which hypotheses are the families $\Cal X_1$ and $\Cal X_2$ locally isomorphic at $0$?
\endproclaim

By {\it locally isomorphic}, we mean that there exists an open neighborhood $W$ of $0\in U$ (respectively in $V$),
and a biholomorphism
$\Phi$ (respectively a CR-isomorphism) between $\Cal X_1(W)=\pi_1^{-1}(W)$ and $\Cal X_2(W)=\pi_2^{-1}(W)$ such that the following
diagram is commutative.
$$
\CD
\Cal X_1(W) \aro >\Phi >> \Cal X_2(W)\\
\aro V \pi_1 VV \aro VV \pi_2 V \\
W \aro >> \text{ Identity }> W
\endCD
$$

We are also interested in the following broader problem.

\proclaim{Question 2} Under which hypotheses are the families $\Cal X_1$ and $\Cal X_2$ locally equivalent at $0$?
\endproclaim

By {\it locally equivalent}, we mean that there exist open neighborhoods $W_1$ and $W_2$ of $0\in U$ (respectively in $V$),
a biholomorphism $\phi$ between $W_1$ and $W_2$ (respectively a diffeomorphism) and a biholomorphism
$\Phi$ (respectively a CR-isomorphism) between $\Cal X_1(W_1)=\pi_1^{-1}(W_1)$ and $\Cal X_2(W_2)=\pi_2^{-1}(W_2)$ such that the following
diagram is commutative.
$$
\CD
\Cal X_1(W_1) \aro >\Phi >> \Cal X_2(W_2)\\
\aro V \pi_1 VV \aro VV \pi_2 V \\
W_1 \aro >> \phi > W_2
\endCD
$$

In other words, $\Cal X_1$ and $\Cal X_2$ are locally equivalent at $0$ if $\phi^*\Cal X_2$ and $\Cal X_2$ are locally isomorphic for
some local biholomorphism $\phi$ of $U$ (respectively $V$) fixing $0$.
\medskip

Fix a family $\Cal X_1$. In this paper, we shall say that this family has the {\it local isomorphism property} (at $0$),
respectively has the {\it local equivalence property} (at $0$) if every other family
$\Cal X_2$ which is pointwise isomorphic to it is locally isomorphic to it (at $0$), respectively locally equivalent (at $0$).
\medskip

It is known since Kodaira-Spencer (see \cite{K-S2}, \cite{We} and Section V.1 of this article) that there exist pointwise isomorphic
families of primary Hopf surfaces which are not locally isomorphic, both in the differentiable and the holomorphic cases. 
\medskip

On the other hand, the classical Fischer-Grauert Theorem \cite{F-G}, can be restated as follows. Let $X$ be a compact complex
manifold and $U$ be an open neighborhood of $0$ in some $\Bbb C^p$. 
Then every trivial family $X\times U$ has the local isomorphism property. This works also for differentiable 
families. Indeed the proof given in \cite{F-G} for holomorphic families is easily adapted to the differentiable case, the core of the
proof being Theorem 6.2 of \cite{K-S1} which is valid both for differentiable and holomorphic families. 
\medskip

Moreover, J. Wehler proved in \cite{We} that, over a smooth base,
 holomorphic families of compact complex tori (in any dimension) as well as holomorphic
families of compact manifolds with negatively curved holomorphic curvature (this implies that they are Kobayashi hyperbolic) 
have the local isomorphism property. This time, the proofs do not adapt to
the differentiable case. 
\medskip
Observe that in the two previous cases, the function $h^0(t)$, that is the dimension of the cohomology group
$H^0(X_t,\Theta_t)$ (where $\Theta_t$ is the sheaf of holomorphic vector fields along $X_t$) is constant for all $t\in U$. It is equal to $n$ in the case of $n$-dimensional tori, and to $0$ in the case of negatively curved manifolds.
Wehler asks in the introduction of \cite{We}
if this condition
is sufficient to have the local isomorphism property.
\medskip

In this paper, we prove that this is the case, even over a singular base. Namely,

\proclaim{Theorem 3}
If $U$ is reduced and if the function $h^0$ is constant for all the fibers of a holomorphic deformation family
$\pi : \Cal X\to U$, then $\Cal X$ has the local isomorphism property.
\endproclaim

We then give examples (both in the differentiable and holomorphic setting) of families not having the local equivalence property, as
well as of locally equivalent but not locally isomorphic families. We classify these counterexamples into two types, 
and we give in Theorem
4 a complete classification of $1$-dimensional families of type II
not having the local equivalence property.
\medskip

Coming back to the search for a criterion, we prove that, surprisingly,
things are completely different in the differentiable case.

\proclaim{Theorem 5} 
There exist differentiable families of 2-dimensional compact complex tori parametrized by an interval that are pointwise isomorphic but not locally isomorphic
at a given point.
\endproclaim

To solve the uniformization problems stated above, we first study the geometry of the Kuranishi space $K$ of a compact complex manifold 
$X$. We show in Theorem 1 
that it has a natural holomorphic foliated structure: two points belonging to the same leaf correspond to biholomorphic complex
structures. More precisely, $K$ admits an analytic stratification such that each piece of the induced decomposition 
(see Section III for more details) is foliated. The leaves are complex manifolds, but the transverse structure of the foliation may
be singular (this happens when the Kuranishi space is singular).
\medskip

The foliation may be of dimension or of codimension zero. In Theorem 2, we prove that there exists leaves of positive dimension (that
is the foliation has positive dimension on some piece of the decomposition) if and only if the function $h^0$ is not constant in the neighbourhood of $0$ in $K$ ($0$ representing the central point $X$). In particular, in many
examples, the foliation is a foliation by points.
\medskip

Although Theorems 3, 4 and 5 on the uniformization problems are not strictly speaking a consequence of Theorems 1 and 2 on the foliated
structure of $K$, the geometric picture of $K$ they bring played an essential role in the understanding and resolution of the
problem. The key ingredients to prove the Theorems are some trivial remarks on diffeomorphisms of the Kuranishi space (see Section II),
the Fischer-Grauert Theorem \cite{F-G}, a result of Namba \cite{Na}  
and a fundamental proposition proved by Kuranishi in \cite{Ku2}.
\medskip
We end the article with a discussion of the relationship between the uniformization problem and the universality of the Kuranishi space.
\vskip1cm
\head
{\bf I. Notations and background}
\endhead
  
Let $X$ be a compact complex manifold. We denote by $X^{\text{diff}}$ the underlying smooth manifold and by $J$ the corresponding
complex operator.
\medskip
A {\it (holomorphic) deformation family} of $X$ is a proper and flat projection $\pi$ from a $\Bbb C$-analytic space $\Cal X$
(possibly non-reduced) over an analytic set $U$ defined on an open neighborhood of $0$ in some $\Bbb C^p$. 
A {\it differentiable deformation family} (see \cite{K-S1}) is a smooth submersion $\pi$ from a smooth
manifold $\Cal X$ endowed with a Levi-flat integrable almost CR-structure over an open
neighborhood $V$ of $0$ in some $\Bbb R^p$, whose level sets
are tangent to the CR-structure.
If the almost CR-structure on $\Cal X$ is not supposed to be integrable, one has a {\it differentiable deformation family of almost-complex structures} of $X$.
\medskip
 In the three cases, the central fiber $X_0=\pi^{-1}(\{0\})$ is assumed to be biholomorphic to $X$.
Sometimes, we consider {\it marked} deformation families of $X$, that is we fix a precise holomorphic identification $i:X\to X_0$. 
\medskip
Let us recall some features of the construction of the Kuranishi space following \cite{Ku1}.
The set of almost complex structures close to $J$ is identified with a
neighbourhood $A$ of $0$ in the space $A^1$ of $(0,1)$-forms on $X$ with values in $T^{1,0}$. In particular, $0$ represents the 
complex structure $J$ we started with (here and in the sequel, the topology used on spaces of sections of a vector bundle over $X$ is induced by some Sobolev norms, see
\cite{Ku1} for more details).
\medskip

Put a hermitian metric $h$ on $X$. Then we have a $\bar\partial$-operator on $A^p$, the space of 
$(0,p)$-forms with values in $T^{1,0}$, a formal adjoint operator $\delta$ with respect to the induced hermitian product on $A^p$ and a Laplace-like operator $\square$.
Let $SH^1$ denote the set of $\delta$-closed forms in $A^1$. Kuranishi proves in \cite{Ku1}

\proclaim{Proposition K1}
For $A$ small enough, there exists a neighborhood $B$ of $0$ in $SH^1$ and  an application $\Xi$ from $A$ to $B$ mapping an almost complex structure $\alpha$ onto a 
$\delta$-closed representant $\Xi(\alpha)$. Moreover, if $\alpha(t)$ is a smooth family of almost complex structures, then so is 
$\Xi(\alpha(t))$.
\endproclaim

By representant, we mean that $\Xi(\alpha)$ and $\alpha$ induce isomorphic almost complex structures on $X^{\text{diff}}$. 
\medskip

Then Kuranishi defines a holomorphic map $G$ from $A^1$ to $A^1$ and proves that it is a
biholomorphism between a special subset of $A^1$ (containing in particular all integrable almost complex structures close
enough to $0$) and
a neighborhood $W$ of $0$ in $H^1$, the (finite-dimensional) space of harmonic forms in $A^1$; or, 
using the Dolbeault isomorphisms, in $H^1(X,\Theta)$. 
\medskip

The {\it Kuranishi space} $K$ is the set of integrable structures in $W$. It is an analytic set. The {\it Kuranishi family} $\Cal K$ 
is the family obtained by endowing each fiber $X^{\text{diff}}\times\{\alpha\}$ of $X^{\text{diff}}\times K\to K$ with the 
corresponding complex structure encoded by $\alpha$. Moreover, it is a marked family.

\medskip
The Kuranishi family is {\it complete} for smooth deformation families as well as for holomorphic deformation families (unmarked and marked), 
not only at $0$ but also at all points of $K$ (shrinking $K$ if necessary). Hence every deformation family $\Cal X$ of $X$ is locally isomorphic to the pull-back of $\Cal K$ by some map $f$ from its base space to the 
Kuranishi space. By abuse of notation, 
we write $\Cal X=f^*\Cal W$. If $\Cal X$ has a marking, then we ask the pull-back to respect the markings.
\medskip

The Kuranishi family is {\it versal} at $0$, i.e. complete in the previous sense with the additional property that
its Zariski tangent space has dimension equal to the dimension of $H^1(X,\Theta)$.
This last property may not be true at points different from $0$.
The versality property is equivalent to the following. Given a holomorphic {\it marked} deformation family $\pi : \Cal X\to U$, 
in the writing $\Cal X=f^*\Cal K$, 
then $f$ may not be unique, but its differential at $0$ is. The same property holds for smooth deformation families.
It must be stressed that this property is related to the markings. It is usually lost when dealing with unmarked families. The marking is necessary
in order to prevent from reparametrizing the family by an automorphism which acts non-trivially in the central fiber.

\remark{Remark}
In the previous setting, if $f$ is unique, then $K$ is called {\it universal}.
The universality property does not hold for any Kuranishi space, see \cite{Wa} and Section V.5.
\endremark

\remark{Remark}
The versality property of $K$ at $0$ implies its unicity (as a germ), see \cite{Ca, Proposition 5.3}.
\endremark

\remark{Remark}
The Kuranishi space may be not reduced at every point \cite{Mu}. This explains why one also considers holomorphic families over a non-reduced base.
In this context, we would like to point out the following subtle point. If $\Cal X$ is a holomorphic family
over a {\it reduced} base, or is a differentiable family, then, in the writing $\Cal X=f^*\Cal K$, the morphism $f$ is in fact a morphism from the base $U$ of $\Cal X$ to
$K^{\text{red}}$, the reduction of $K$. Hence the fact that $K$ is reduced or not is not relevant for these families. 
But, if $U$ is non-reduced, then such a $f$ is not completely determined by its image in $K^{\text{red}}$; one has also to specify the value of its differential
at each point. Anyway, keeping in mind this difference, it is still true that such a family is obtained by pull-back from the Kuranishi family and that the differential of the pull-back map at $0$ is unique; that is,
the notions of completeness and versality remain the same. 
In particular, 
from the previous remark, one deduces that if a deformation family of $X$ over some reduced analytic space $A$ is versal for differentiable families, then it is diffeomorphic
to the reduction of the Kuranishi space of $X$.
\endremark
\medskip
 
Finally, we will make an intensive use of the following Proposition

\proclaim{Proposition K2} 
If $h^0$ is constant on $K$, there exists a neighborhood $U$ of $0$ in $K$ and a neighborhood
$\Cal U$ of the identity in the group of diffeomorphisms of $X^{\text{diff}}$ such that, for all couples
$(t,t')$ of distinct points of $U$, we have

$$
\alpha_{t'}\equiv f_* \alpha_t \quad\text{ for some diffeomorphism } f \Longrightarrow f\not \in \Cal U
$$
\endproclaim

Here we have $\alpha_t=G^{-1}(t)$ (respectively $\alpha_{t'}=G^{-1}(t')$).
To fully understand this statement, recall that the Kuranishi family is constructed as a families of complex operators. Hence every
fiber $K_t$ is naturally defined as $(X^{\text{diff}},\alpha_t)$. This crucial Proposition is proved by Kuranishi in \cite{Ku2} and used to show that $h^0$ constant implies the universality of $K$.
We will discuss this at the end of the article.

\vskip1cm
\head
{\bf II. Preliminary remarks on diffeomorphisms of the Kuranishi family}
\endhead

Let us begin with some definitions.

\definition{Definition} 
A {\it diffeomorphism} of the Kuranishi space $K$ is a bijective map $\phi$ from some open neighborhood of $0$ in its reduction $K^{\text{red}}$ onto 
some open neighborhood of $0$ in $K^{\text{red}}$ such that 
\medskip
\noindent (i) It sends a complex structure onto an isomorphic complex structure.

\noindent (ii) Both $\phi$
and $\phi^{-1}$ are restrictions to $K^{\text{red}}\cap W'\subset W'\subset W\subset H^1(X,\Theta)$ of a 
smooth map of $W'\subset W$ for some $W'$.
\enddefinition
Such a diffeomorphism is generally, but not always, assumed to fix $0$.
Notice that such a map is smooth in the sense of \cite{Ku1}.

\definition{Definition}
A {\it diffeomorphism} of the Kuranishi family $\Cal K$ is a continuous map $F$ from some open neighborhood of $K_0$ in $\Cal K$
to $\Cal K$ such that

\medskip

\noindent (i) $F$ descends as a diffeomorphism $f$ of $K$.

\noindent (ii) the restriction of $F$ to any fiber of $\Cal K\to K$ is a biholomorphism.

\noindent (iii) $F$ is CR in the following sense. Since $\Cal K\to K$ is a flat morphism, it is locally isomorphic at each point
to an open set of $\Bbb D^{\dim K_0}\times K$. Representing $F$ locally as a map between two such charts, we ask it to be
holomorphic in the $\Bbb D^{\dim K_0}$-variables, and smooth in the other variables (in the sense of the previous definition).
\enddefinition

Notice that, even when $F$ fixes the central fiber, {\it we do not ask it to respect the markings}.
\medskip
In the same way, we define {\it automorphisms} of $K$ as isomorphisms of some open neighborhood of $0$ in the analytic space $K$ (and not $K^{\text{red}}$ this time)
generally 
fixing $0$ (and thus as restrictions to $K$ of local isomorphisms of $W$ at $0$); and {\it automorphisms} of $\Cal K$ as local
isomorphisms of $\Cal K$ descending as automorphisms of $K$. Finally, all these definitions apply with trivial changes
to other deformation families of $X$ than the Kuranishi family.

\remark{Remark}
In the definition of a diffeomorphism (and an automorphism) of $K$, we consider $K$ as an analytic space of a $\Bbb C$-vector space,
and not as a set of complex operators. This explains why a diffeomorphism of $K$ may not lift to a diffeomorphism of $\Cal K$.
This lifting problem is very close to the local isomorphism problem. Indeed, we will give a criterion for lifting an automorphism
of $K$ in Corollary 4, as a consequence of the criterion to have the local isomorphism property. And we also give in Lemma 6 an
example of an automorphism of $K$ which does not lift.
\endremark
\medskip

In the second part of this Section, we deal with the problem of extending a diffeomorphism (respectively an automorphism)
of the central fiber $K_0$ to a diffeomorphism (respectively an automorphism) of $\Cal K$. Let us make first the following
trivial remark. Let $\phi$ be an automorphism of $X$. Via the marking of $\Cal K$, we consider it as an automorphism of
$K_0$. The family $(\Cal K, K)$ with the new identification $\phi \circ i$ is a new versal
family for $X$, hence, by unicity, there exists an automorphism $\Phi$ of $\Cal K$ fixing $0$ and extending $\phi$. 
\medskip

The two following Lemmas are trivial but of fundamental importance for the sequel. Part (i) of the first one is even weaker than the
previous statement but it has the advantage to admit slight generalizations stated in Lemma 2 and Proposition 1. 

\proclaim{Lemma 1}

(i) Let $\phi$ be an automorphism of $X$. Then there exists a diffeomorphism $\Phi$ of $\Cal K$ extending $\phi$.

(ii) Let $\phi$ be a diffeomorphism of $X^{\text{diff}}$ such that $\phi_* 0$ belongs to the set $A$ of Proposition K1. Then there exists 
a diffeomorphism $\Phi$ of $\Cal K$ extending $\phi$.
\endproclaim

\demo{Proof}
In both cases, see $\phi$ as a diffeomorphism of $X^{\text{diff}}$. Then it satisfies $\phi_*0= 0$ (case (i)) or $\phi_* 0$ close
to $0$ (case (ii)).  In other words, this diffeomorphism induces a map
$$
\phi_*\ :\ \alpha\in A^1\longmapsto \phi_*\alpha\in A^1
$$
with $0$ as fixed point or with $0$
sent close to $0$. Consider the following composition
$$
W'\subset W\subset H^1\aro> G^{-1}>> A^1\aro >\phi_* >> A^1\aro>\Xi >>SH^1\aro >G>> H^1
$$
taking $W'$ small enough to have $\phi_*(W')\subset A$.
\medskip

This gives $\tilde\phi$, a map from $W'\subset W$ to $H^1$. This $\tilde\phi$ respects the almost complex structures, that is sends 
an almost complex structure onto one which is isomorphic. 
Hence, it induces naturally a smooth map from $\Cal K$ to $\Cal K$ that we denote by $\Phi$.
\medskip

Consider the image $K'\subset K^{\text{red}}$ of $\Phi$. Since $\Phi$ respects the complex structures, we have that $\Cal K^{\text{red}}$ is diffeomorphic to 
$\Phi^*\Cal K'=\Phi^*\Cal K^{\text{red}}$. Hence, the analytic set $K'$ is versal for differentiable families, so, as remarked at the end of Section I, $K'$ and $K^{\text{red}}$ must be equal (as germs of $0$) and $\Phi$
must be a diffeomorphism.
 $\square$
\enddemo

\proclaim{Lemma 2}
Let $\phi$ be an automorphism of $X$ isotopic to the identity. Then there exists a diffeomorphism $\Phi$ of $\Cal K$ extending $\phi$
isotopic to the identity. Moreover, if we fix an isotopy $\phi_t$ between $\phi$ and the identity on $X$, then 
the isotopy between $\Phi$ and the identity of $\Cal K$ can be assumed to be equal to $\phi_t$ when restricted to $X$.
\endproclaim

\demo{Proof}
Apply the proof of Lemma 1 to each member of the isotopy $\phi_t$. We thus obtain a family $\Phi_t$ of diffeomorphisms
of $\Cal K$ extending $\phi_t$ for all $t$. By compacity of $[0,1]$, all the $\Phi_t$ can be defined on a same open neighborhood
of $0$. Finally smoothness in $t$ comes from Proposition K1.
$\square$
\enddemo

Of course, the important point in Lemma 2 is that $\Phi$ is isotopic to the identity.
We draw from these lemmas an important consequence.

\proclaim{Proposition 1}
Assume that the function $h^0$ is not constant at $0\in K$. Then, we can find
$\phi$, an automorphism of $X$ isotopic to the identity, such that 
\medskip

\noindent (i) it extends as a diffeomorphism $\Phi$ of $\Cal K$ isotopic to the identity.

\noindent (ii) the projection of every such extension on $K$ gives a diffeomorphism of $K$ whose germ at $0$ is not the germ
of the identity.
\endproclaim

\demo{Proof}
The first step of the proof is a very classical argument. Assume without loss of generality that $K$ is reduced.
The function 
$h^0$ is known to be upper semi-continuous. The assumption that it is not constant in a neighborhood of $0$ implies then that it has
a strict maximum at $0$. Take a basis of $H^0(K_0,\Theta_0)$. If every vector field of this basis could be extended to a vector field
of the family $\Cal K$ (tangent to the fibers, globally smooth and holomorphic along the fibers), 
at a point $t\in K$ close enough to $0$, they all would form a free family of dimension $h^0(0)$ of
$H^0(K_t,\Theta_t)$, contradicting the inequality $h^0(t)<h^0(0)$.
\medskip
Hence, there exists a global vector field $\xi$ 
on $K_0$ which cannot be extended as a vector field of $\Cal K$ tangent to the fibers. Let
$\phi$ be the corresponding automorphism of $K_0$ isotopic to the identity obtained by exponentiation for small time.
By Lemma 2, there exists a diffeomorphism $\Phi$ of $\Cal K$ extending $\phi$
isotopic to the identity, proving (i).
\medskip
Now, for every such choice of $\Phi$, the induced diffeomorphism of $K$ cannot be the identity, even in germ, otherwise the global vector field
of $\Cal K$ obtained by differentiating $\Phi$ would be tangent to the fibers and would extend $\xi$. Contradiction which
proves the Proposition.
$\square$
\enddemo

A classical result of \cite{K-S1} states that, if the function $h^0$ is constant in a smooth deformation family $\pi : \Cal X\to V$,
then every automorphism  of $X_0$ isotopic to the identity can be extended as a diffeomorphism of the family $\Cal X$ which is the
identity on $V$. Automorphisms of $X$ which does not extend as automorphisms of $\Cal K$ that are the identity on $K$ are usually called 
{\it obstructed automorphisms}. Proposition 1 tells us that obstructed automorphisms extend as diffeomorphisms of $\Cal K$ with
non-trivial projection on $K$, but isotopic to the identity.

\vskip1cm

\head
{\bf III. Foliated structure of the Kuranishi space}
\endhead

\vskip .5cm

\subhead
1. Local submersions
\endsubhead

Let $t$ be a point of $K$ corresponding to the complex manifold $X_t=(X^{\text{diff}}, J_t)$. 
Denote by $(K)_t$ the space $K$ but with base point $t$ and not $0$ and by $(\Cal K)_t$ the
corresponding deformation family of $X_t$ (choosing some identification maps). The family $(\Cal K)_t\to (K)_t$ is complete at $t$, but not always versal. On the
other hand, let $K(t)$ be the Kuranishi space of $X_t$, and $\Cal K(t)$ the corresponding versal family. We thus have a sequence of
pointed analytic spaces
$$
(K(t),0) \aro >i_t >> ((K)_t, t) \aro > s_t >> (K(t),0)
\leqno (\Cal S)
$$
which lifts to a sequence of maps between families
$$
\Cal K(t) \aro > I_t >> (\Cal K)_t \aro > S_t >> \Cal K(t)
$$

And $(\Cal S)$ is the restriction of the sequence (defined on neighborhoods of the base points)
$$
(H^1(X_t,\Theta_t),0) \aro > \tilde \imath_t >> (H^1(X_0,\Theta_0), \alpha_t) \aro > \tilde s_t >> (H^1(X_t,\Theta_t), 0)
$$

\proclaim{Lemma 3}
The map $\tilde s_t$ is a submersion at $\alpha_t$.
\endproclaim

\demo{Proof}
Since $K(t)$ is versal at $0$, the composition 
$\imath_t \circ s_t$ is a local isomorphism at $0$. So is $\tilde \imath_t \circ \tilde s_t$.
Hence $\tilde s_t$ is a submersion at $J_t$.
$\square$
\enddemo

Apply the submersion Theorem to $\tilde s_t$. This gives a diagram
$$
\CD V\subset (H^1(X_0,\Theta_0), J_t) \aro > \tilde s_t >> W\subset (H^1(X_t,\Theta_t), 0)\\
\aro V \text{local biholomorphism} VV \aro AA \text{natural projection} A \\
W\times B \aro >> \text{identity} > W\times B \\
\endCD
$$
where $B$ is the unit euclidean ball of $\Bbb C^p$ for $p=h^1(0)-h^1(t)$ and $h^1(t)$ denotes the dimension of
$H^1(X_t,\Theta_t)$.

\medskip

This submersion allows to locally foliate $H^1(X_0,\Theta_0)$ in a neighborhood of $\alpha_t$. The leaves correspond to 
deformation families of almost complex structures which are pull-back by $\tilde s_t$ of a constant family. In other words,
the points of a same leaf all define the same almost complex structure up to isomorphism.
\medskip

When we restrict to $K$, we obtain the diagram of analytic spaces
$$
\CD (K)_t \aro > s_t >> K(t)\\
\aro V \text{local biholomorphism} VV \aro AA \text{natural projection} A \\
K(t)\times B \aro >> \text{identity} > K(t)\times B \\
\endCD
$$
that defines a local foliated structure of $K$.
\medskip
We aim at gluing those local foliations together into a global foliation. This brings some problems since the induced foliations
in two arbitrarily close points may be of different dimensions. To overcome this problem, it is necessary to decompose the space $K$.
\vskip .5cm

\subhead
2. Decomposition of $K$
\endsubhead

Decompose $K$ as a disjoint union
$$
K=K^{\text{min}}\sqcup\hdots\sqcup K^{\text{max}}
\leqno (\Cal D)
$$
where $K^i$ denotes the set of points $t$ of $K$ such that the dimension of $K(t)$ is $i$. Observe that the completeness of
$K$ {\it at each point} implies that $0$ belongs to $K^{\text{max}}$. Set $Z_j=\sqcup _{i\leq j} K^i$. The sequence $Z_{\text{min}}\subset \hdots\subset Z_{\text{max}}=K$ is a stratification of $K$.
\medskip
We want to show that this decomposition is analytic in the sense that,  for all $\text{min}\leq i<\text{max}$, the set $Z_i$ is an analytic
subset of $Z_{i+1}$ and thus of $K$.

\proclaim{Lemma 4}
The decomposition $(\Cal D)$ of $K$ is analytic.
\endproclaim

\demo{Proof}
Assume first that $K$ is reduced. Define
$$
E_c=\{t\in K \quad\vert\quad h^1(t)\geq c\}
\leqno \text{min}\leq c\leq \text{max}
$$
These sets forms an analytic stratification of $K$, see \cite{Gr}. Call $\Cal G=(G_c)$ the associated decomposition.
\medskip
On the other hand, denoting by $\Cal G_c$ the family of complex structures with base $G_c$, the results of \cite{Gr}
show that the group $H^1(\Cal G_c,\Theta)$ is isomorphic to a locally free sheaf over $G_c$ whose stalk at $t$ is $H^1(X_t,\Theta_t)$.
The set of integrable structures in this sheaf is given as the zero set $Z_t$ of a field of quadratic forms in the fibers which is
analytic in $t$.
\medskip
This allows us to analytically stratify each piece $G_c$ by 
$$
G_{c,d}=\{t\in G_c\quad\vert\quad \text{codim } Z_t\leq d\}
\leqno d\leq c
$$

Then the $E_c$ respectively the $G_{c,d}$ are analytic sets of $K$ respectively of $G_c$, 
whereas $G_c\subset E_c$ is a quasi-analytic open set. 
\medskip

Observe now that the completeness property of a Kuranishi space {\it at each point} close enough from the base point implies
that the function
$t\in K\longmapsto \dim K(t)$
is upper semi-continuous for the standard topology.  Since $G_{c,d}$ is an analytic set of $G_c$ and since the Zariski
closure of $G_c$, that is $E_c$, is equal to its closure for the standard topology, this proves that the Zariski closure 
$\overline{G_{c,d}}$
of each $G_{c,d}$ in $\overline{G_c}=E_c$ is just
$$
\overline{G_{c,d}}=\{t\in E_c\quad\vert\quad \text{codim } Z_t\leq d\}
$$

We now just have to set
$$
F_i=\cup_{c-d\geq i} \ \overline{G_{c,d}}
\leqno \text{min}\leq i\leq \text{max}
$$
to obtain an analytic stratification of $K$ whose associated decomposition $F_{i}\setminus F_{i+1}$ ($i>\text{min}$) and $F_{\text{min}}$
is the decomposition $\Cal D$.
\medskip
If $K$ is not reduced, we just perform the previous stratification on its reduction, then we obtain the decomposition
$(\Cal D)$ by putting on each $F_i$ the multiplicity induced from $K$. Of course, we may be forced to add more pieces if some component contains an analytic
subspace of higher multiplicity (think of a double point inside a line; then $K$ has two components, whereas $K^{\text{red}}$ has just one).
$\square$
\enddemo

Let $t\in K^i$ for some $i$. By definition, $i_t$ and $s_t$ respect the decompositions of $K$ and $K(t)$, that is
$$
i_t(K(t)^{\text{max}})=K^i \qquad\text{et}\qquad s_t(K^i)=K(t)^{\text{max}}
$$

It is thus possible to restrict the submersions $s_t$ to a piece of the decomposition to obtain the following diagram

$$
\CD (K^i)_t \aro > s_t >> K(t)^{\text{max}}\\
\aro V \text{local biholomorphism} VV \aro AA \text{natural projection} A \\
K(t)^{\text{max}}\times B \aro >> \text{identity} > K^{\text{max}}_t\times B \\
\endCD
$$ 

In the sequel, $s_t$ will always denote the restricted submersion from $K^i$ to the piece $K(t)^{\text{max}}$, for $t\in K^i$.

\vskip .5cm
\subhead
3. Foliated structure of $K$
\endsubhead

Using the submersion $s_t$, we will define a foliated structure on each piece $K^i$.

\definition{Definition}
Let $X$ be an analytic space. A {\it transversally singular foliation} of dimension $p$ on $X$ is given by an open covering
$(U_\alpha)$ of $X^{\text{red}}$ and local isomorphisms $\phi_\alpha : U_\alpha\to B\times Z_\alpha$ (for $B$ the unit ball of $\Bbb C^p$ and
$Z_\alpha$ a reduced analytic space) such that the changes of charts $\phi_{\alpha\beta}\equiv \phi_\beta\circ (\phi_\alpha)^{-1}$
preserve the plaques $B\times\{pt\}$.
\enddefinition

Thus a transversally singular foliation is a lamination with a transverse structure of an analytic space. The choice for this somewhat 
unusual name (rather than analytic lamination or something analogous) comes from the fact that we find more judicious to reserve
the word lamination to a situation where the total space has no analytic structure.
\medskip 
Now, given such a foliation, one may define the leaves as in the classical case by gluing the plaques. Hence the leaves are holomorphic
manifolds.
Remark also that the germ of the analytic space $Z_\alpha$ is the same along a fixed leaf.
\medskip
Starting from a non-reduced space $X$, one may consider such a foliation on it as a foliation of its reduction with holomorphic submanifolds as leaves; or one may endow
$X$ itself with a "non-reduced" foliation whose leaves are non-reduced holomorphic submanifolds. Both points of view are equivalent. 
 \medskip
We may state

\proclaim{Theorem 1} Let $K$ be the Kuranishi space of a manifold $X$. Consider the decomposition $\Cal D$ of $K$.

Then each piece $K^i$ admits a transversally singular foliation $\Cal F^i$ locally defined by the submersions $s_t$ of Section 2.
\endproclaim

Notice that two points belonging to the same leaf correspond to the same complex manifold (up to biholomorphism).
Notice also that the Kuranishi family $\Cal K$ admits an induced decomposition in pieces $\Cal K^i$, each of these pieces being
foliated. By \cite{F-G}, the leaf of $\Cal K^i$ corresponding to $t\in K^i$ is a locally trivial fibre bundle with fibre $X_t$
over the leaf of $K^i$
through $t$. And the foliation of $K$ (respectively $\Cal K$) extends to a holomorphic foliation, respectively an almost-complex
foliation (this time in the classical sense,
that is with smooth transverse structure) of $W$ and respectively $\Cal W$ using the submersions $S_t$.
To simplify the exposition, we will always consider the foliation of $K$.

\demo{Proof}
Take a covering of a piece $K^i$ by open sets where a submersion $s_t$ is well defined. Using the submersion Theorem, we obtain foliated charts modelled on a product of a ball of dimension $(\dim K -i)$ with an analytic space of dimension $i$. Now, the changes of charts
respect the leaves, since they have an intrinsic geometric definition: the leaf through $t\in K^i$ is the maximal connected subset of
$K^i$ of points corresponding to the complex manifold $X_t$ up to biholomorphism.
$\square$
\enddemo

Remark that this proof adapts immediately to the case of $W$, with the only difference that $X_t$ may just be an almost complex manifold.
\medskip

\definition{Notation and Definition}
We denote by $\Cal F$ the global foliation of Theorem 1, that is the union of the foliations $\Cal F^i$.

We say that the foliation $\Cal F$ is {\it trivial} if each foliation $\Cal F^i$ is a foliation by points. 
\enddefinition

The foliation $\Cal F$ is trivial if and only if the decomposition $(\Cal D)$ has a single piece, that is if and only if the dimension
of $K(t)$ is constant  near $0$. On the other hand, observe that $\Cal F$ may be of codimension $0$ (that is each $\Cal F^i$ has 
codimension $0$ in $K^i$) with a non-trivial decomposition (see the examples below).

\vskip .5cm
\subhead
4. Examples
\endsubhead

\noindent (i) Let $X$ be a complex torus of dimension $n$. Following \cite{K-S2}, the Kuranishi space of $X$ may be represented by
a neighborhood of $0$ in $\Bbb C^{n^2}$ and it is versal at each point. Therefore, the decomposition $(\Cal D)$ has a unique piece
and  the foliation $\Cal F$ is trivial.
\medskip

However, observe that given a well-chosen compact complex torus, there exists an infinite sequence of points of $K$ corresponding to 
this torus \cite{K-S2, p. 413}. This shows that two points belonging to different leaves of $\Cal F$ may nevertheless define
the same complex structure.

\bigskip
\noindent (ii) Let $X$ be the Hirzebruch surface $\Bbb F_2$. Its Kuranishi space may be represented by a unit $1$-dimensional disk
whose non-zero points all correspond to $\Bbb P^1\times\Bbb P^1$, see \cite{Ca}. The decomposition is
$$
K=K^1\sqcup K^0=\{0\}\sqcup \Bbb D^*
$$
and both foliations have codimension zero.
\bigskip

\noindent (iii) Let $X$ be the Hopf surface obtained from $\Bbb C^2\setminus \{(0,0)\}$ by taking the quotient by the group 
$\langle 2\text{Id} \rangle$ generated
by the homothety $(z,w)\mapsto 2\cdot (z,w)$. Its Kuranishi space is described in \cite{K-S2}. It may be represented by a neighborhood
of the matrix $2\text{Id}$ in
$$
K=\{A\in \text{GL}_2(\Bbb C) \quad\vert\quad \vert \text{Tr } A \vert >3,\ \vert \Delta (A)\vert=\vert (\text{Tr } A)^2- 4\det A\vert<1
 \}
$$

A point $A$ of $K$ corresponds to the Hopf surface $\Bbb C^2\setminus \{(0,0)\}/\langle A\rangle$.
If $A$ is a multiple of the identity, then the corresponding Kuranishi space $K(A)$ has dimension four; in other words $K$ is versal
along the set
$$
\Delta=\{\lambda \text{Id} \quad\vert\quad \vert\lambda\vert >\dfrac{3}{2}\}
$$

However, if $A$ is not a multiple of the identity, the dimension of $K(A)$ drops to $2$. Thus we decompose $K$ into two pieces
$$
K^4=\Delta \quad\text{ and }\quad K^2=K\setminus\Delta
$$
On the other hand, consider the map
$$
\phi\ :\ A\in K\longmapsto (\text{Tr } A, \Delta(A)) \in\Bbb C^2
$$
Let $(\sigma,\delta)$ be a point of $\Bbb C^2$ with $\vert \sigma\vert >3$ and $\vert \delta\vert <1$. If $\delta$ is different from
zero, all points of $\phi^{-1}(\sigma,\delta)$ correspond to the same Hopf surface. If $\delta$ is zero, the same is true
for all points of $\phi^{-1}(\sigma,\delta)$ except for $\sigma/2\cdot \text{Id}$, which corresponds to a different Hopf surface. Notice that in 
this case, the level set $\phi^{-1}(\{(\sigma,\delta)\})$ is singular at $\sigma/2\cdot \text{Id}$.
\medskip
As a consequence of all that, the foliation $\Cal F^4$ is a foliation by points, whereas the foliation $\Cal F^2$ is a non-trivial 
one, which is given by the level sets of the submersion $\phi$ restricted to $K^2$. It has dimension and codimension two.

\vskip 1cm

\head
{\bf IV. Non-triviality criterion for $\Cal F$}
\endhead

The aim of this section is to prove the following result.

\proclaim{Theorem 2} 
Let $K$ be the Kuranishi space of $X$ and let $\Cal F$ be the foliation of $K$ constructed in section III.

Then $\Cal F$ is trivial if and only if $h^0$ is a constant function on $K$.
\endproclaim

This leads to the following corollary.

\proclaim{Corollary 1}
The Kuranishi space $K$ is versal at all points if and only $h^0$ is a constant function on $K$.
\endproclaim

\demo{Proof of Corollary 1}
Combine Theorem 2 and the remark after the definition of triviality for $\Cal F$.
$\square$
\enddemo

Let us proceed to the proof of Theorem 2.

\demo{Proof of Theorem 2}
Assume $h^0$ constant on $K$ and assume at the same time that $\Cal F$ is non-trivial. Thus there exists a piece $K^i\subset K$
whose foliation $\Cal F^i$ has positive-dimensional
leaves. So there exist non-constant smooth paths $c : [0,1]\to K$ such that the induced family $\Cal C=c^*\Cal K$ has all fibers
biholomorphic. Choose such a non-constant path staying inside the neighborhood $U$ appearing in 
Proposition K2. Now Fischer-Grauert Theorem \cite{F-G} implies that $\Cal C$ is the trivial family; in other words there exists
$(\phi_t)_{t\in [0,1]}$ an isotopy such that 
\medskip
\noindent (i) $\phi_0\equiv \text{Id}$.

\noindent (ii) For all $t\in [0,1]$, we have $(\phi_t)_*\alpha_{c(0)}=\alpha_{c(t)}$.
\medskip

For $t$ small enough, we have $\phi_t$ in $\Cal U$, violating Proposition K2. Contradiction. The foliation is trivial.
\medskip
Reciprocally, assume $h^0$ non-constant. Then by Proposition 1, there exists an automorphism $\phi$ of $X$ isotopic to the identity
all of whose extensions as a diffeomorphism of $\Cal K$ does not project onto the identity of $K$ on any neighborhood of $0$. 
Let $\Phi$ be one of these 
extensions; still by Proposition 1, recall that $\Phi$ may be chosen isotopic to the identity. Let $\Phi_t$ be the isotopy.
All that means that there exist points $x\in K$ arbitrary close to $0$ such that the path $t\in [0,1]\mapsto \Phi_t(x)\in K$
is a non-constant path. But Fischer-Grauert Theorem may be geometrically reformulated as follows.

\proclaim{Lemma 5}
 The Kuranishi space of a compact complex
manifold does not contain a non-constant path passing through $0$ all of whose points correspond to $X$.
\endproclaim

\demo{Proof}
Assume the contrary and consider the family associated to such a non-constant path. By \cite{F-G}, its Kodaira-Spencer map is zero at every point. On the other hand, $K$ is versal at each point of
the path. This is due to the fact that at any such point $t$, the Zariski dimension of $K$ is greater than $h^1(t)$ by completeness. Since $h^1(t)$ is equal to $h^1(0)$, they must be equal yielding the versality 
of $K$ at $t$. Now all that means that this non-constant path should be parametrized by a map whose derivative at each point is zero. Contradiction.
$\square$
\enddemo

As a consequence, for such a point $x$, the space $(K)_x$ cannot be the Kuranishi space $K(x)$. That is it is not versal at $x$.
But we already noticed in III.3 that this is enough to prove that the foliation is non-trivial. 
$\square$ 
\enddemo

Finally, note that:

\proclaim{Corollary 2}
The stratum $K^{\text{max}}$ (or, in the non-reduced case, the union of the strata corresponding to the stratum $K^{\text{max}}$ of $K^{\text{red}}$) is the set of points $t\in K$ 
such that $h^0(t)=h^0(0)$.
\endproclaim

\demo{Proof}
Assume that $K$ is reduced. Let
$$
H=\{t\in K\quad\vert\quad h^0(t)=h^0(0)\}
$$

This is an analytic space by \cite{Gr} (recall that $h^0(t)\geq h^0(0)$ implies equality). Arguing exactly as in the first part of the
proof
of Theorem 2, we show that $\Cal F$ is trivial on $H$. So $H$ is included in $K^{\text{max}}$.
\medskip
Conversely, assume that $h^0$ is not constant on $K^{\text{max}}$. Then arguing as in the second part of the proof of Theorem 2, we show that
$\Cal F^{\text{max}}$ is non trivial. Contradiction.
\medskip
The non-reduced case can be treated in a similar way.
$\square$
\enddemo

\remark{Remark}
Indeed, although Proposition K2 is stated for the complete Kuranishi space $K$, a quick look at the proof shows that it is valid
in restriction to any subset $V\subset K$ where $h^0$ is constant equal to $h^0(0)$.
\endremark

\vskip1cm

\head
{\bf V. The isomorphism and equivalence problems}
\endhead

We refer to the introduction for the definition of these two problems. Let us give two more definitions.

\definition{Definitions}
Let $\Cal X_1$ and $\Cal X_2$ be two families which are pointwise isomorphic but not locally isomorphic at $0$. Then we say that they form
{\it a type (II) counterexample (to the isomorphism property)} if there exist
$$
f,\ g\ :\ (U,0)\longrightarrow (K,0)
$$
such that
\medskip

\noindent (i) We have $\Cal X_1=f^*\Cal K$ and $\Cal X_2=g^*\Cal K$.

\noindent (ii) There exist $U_1\subset U$ and $U_2\subset U$ such that $f(U_1)$ and $g(U_2)$ are equal. 
\medskip
 
And we say that they form a {\it type (I) counterexample} if we cannot find $f$ and $g$ as above.
\enddefinition
Same definitions are valid for the equivalence problem.
Roughly speaking, a type (II) counterexample is a counterexample obtained by reparametrization, whereas a type (I) counterexample
relies on the particular geometric structure of the Kuranishi space.
\vskip .5cm
\subhead
1. Counterexamples
\endsubhead

A counterexample of type (II) for the two problems (in both differentiable and holomorphic cases)
can be found in \cite{K-S2} (and explained in \cite{We}). Start with the Hopf surface $\Bbb C^2\setminus
\{(0,0)\}/\langle 2\text{Id}\rangle$. We use the notations of Section III.4.(iii). Define $\Cal X_1$ as the family corresponding to an
embedding disk (respectively an interval) in the closure of the two-dimensional leaf $\phi^{-1}(4,0)$. This is an example of a 
jumping family. Let us give a precise definition.

\definition{Definition}
A holomorphic (respectively differentiable) {\it jumping family} is a family $\pi :\Cal X\to U$ (respectively $\pi : \Cal X\to V$)
such that
\bigskip

\noindent (i) it is trivial outside $0$, but the central fiber $X_0$ is not biholomorphic
to the generic fiber.

\noindent (ii)  The Kodaira-Spencer map at $0$ is not zero. 
\enddefinition

Notice that the Kodaira-Spencer map $\rho_1$ of our family $\Cal X_1$ is not zero at $0$, since it is an embedding at $0$. This follows from the versality 
property of $K$ at $0$.
Consider now the
ramified covering
$$
z\in\Bbb D\longmapsto z^2\in\Bbb D \quad\text{or more generally }z^n\in \Bbb D
$$
or respectively $t\in I\longmapsto t^3\in I$. Then define $\Cal X_2$ as the pull-back of $\Cal X_1$ by this application. By the
``chain-rule for the Kodaira-Spencer map'', we have at $0$
$$
\displaylines{
\rho_2\left (\dfrac{\partial}{\partial z}\right )=\rho_1 \left (\text{Jac }_0(z\mapsto z^n)\cdot \dfrac{\partial}{\partial z}
\right )=\rho_1(0)=0\hfill\cr
\hfill \text{ respectively }
\rho_2\left (\dfrac{\partial}{\partial z}\right )=\rho_1\left (\text{Jac }_0(t\mapsto t^3)\cdot \dfrac{\partial}{\partial z}
\right )=\rho_1(0)=0 }
$$
so $\Cal X_2$ is not isomorphic, nor equivalent, to $\Cal X_1$.

\medskip
Of course, this construction can be generalized starting from any jumping family (for example, one can take the jumping family with
the Hirzebruch surface $\Bbb F_2$ as central fiber and $\Bbb P^1\times\Bbb P^1$ as generic fiber; this shows that such counterexamples
exist even for projective manifolds). So we state:

\proclaim{Proposition 2}
A (holomorphic or differentiable) jumping family has neither the local isomorphism, nor the local equivalence property.
\endproclaim

It is important for the sequel to observe that the function $h^0$ is not constant in a jumping family (cf \cite{Gri}).
\medskip
We give now a type (I) counterexample for the two problems, in both differentiable and holomorphic cases. 
Although it can be obtained easily from the treatment of Hopf surfaces in \cite{K-S2},
we do not know of any reference where it is described.
\medskip
Once again, we use the results and the notations of Section III.4.(iii). Consider
\bigskip
$$\left\{
\eqalign{
\Cal X_1&=\Bbb C^2\setminus\{(0,0)\}\times\Bbb D/
{\left\langle \pmatrix 
2+t &t \\
0 &2+t
\endpmatrix
,\ t\right\rangle} \cr
\null\cr
\Cal X_2&=
\Bbb C^2\setminus\{(0,0)\}\times\Bbb D/
{\left\langle \pmatrix 
2+t &t^3 \\
0 &2+t
\endpmatrix
,\ t\right\rangle }
}
\right .
$$
\bigskip
Replacing $\Bbb D$ by $I$ in the definition of $\Cal X_1$ and $\Cal X_2$, one obtains a differentiable counterexample.

\medskip

We claim that $\Cal X_1$ and $\Cal X_2$ are pointwise isomorphic but not locally equivalent at $2\text{Id}$, and finally that they have
distinct image in $K$. The last point is a direct consequence of the fact that, since the families are embedded, 
same image would imply locally isomorphic.

\medskip

Now, an elementary computation shows that for $t\not = 0$, the fibers $(X_1)_t$ and $(X_2)_t$ are biholomorphic and conjugated by
$$
P(t)=\pmatrix
\pm t &q \\
0 &\pm t^{-1}
\endpmatrix
$$
where $q$ is any complex number (we assume without loss of generality that $P$ has determinant one). 
Since $(X_1)_0=(X_2)_0$ and since none of these conjugating matrices extend at $0$, we are done for
the isomorphism problem. Finally, for the equivalence problem, just observe that if $t$ and $t'$ are distinct and both different 
from $0$, then $(X_1)_t$ and $(X_2)_t$ are not biholomorphic (look at the traces). Hence, in this case, there is no difference between the isomorphism
problem and the equivalence problem.
\medskip
Notice that $h^0$ is not constant along these families, dropping from $4$ (at $0$) to $2$.
\medskip
Let us give now examples of locally equivalent but not locally isomorphic families. We still use the Kuranishi space of the Hopf 
surface described in Section III, 4, (iii). The key point is given by the following Lemma.

\proclaim{Lemma 6}
The map $A\in K\to \ ^t\kern -.1ptA\in K$ is an automorphism of $K$ fixing $2\text{Id}$ which does not lift to an automorphism of $\Cal K$.
\endproclaim

\demo{Proof}
This is clear for $K$ using the fact that $A$ and $^t\kern -.1pt A$ are conjugated. On the other hand, assume that this automorphism lifts
to an automorphism of $\Cal K$. Then, in a neighborhood of $2\text{Id}$, it would be possible to find a family of invertibles
matrices $P(A)$ depending holomorphically on $A$ such that
$$
^t\kern -.1pt A=P^{-1}(A)\cdot A\cdot P(A)
$$
where we assume without loss of generality that $P(A)$ has determinant equal to one.
Straightforward computations show that we must have
$$
P(A)=\pmatrix
\alpha &\pm i \\
\pm i &0
\endpmatrix
\quad\text{ for }\quad
A=\pmatrix
2 &t \\
0 &2
\endpmatrix
\eqno t\in\Bbb C
$$
where $\alpha$ is any complex number and can be chosen independently of $t$. And we must also have
$$
P(A)=\pmatrix
\alpha &0 \\
  0  &1/\alpha
\endpmatrix
\quad\text{ for }\quad
A=\pmatrix
2 &0 \\
0 &2+t
\endpmatrix
\eqno t\in\Bbb C
$$
where $\alpha$ is any non-zero complex number and can be chosen independently of $t$.
\medskip
Since these two families do not have any common limit where $t$ goes to zero, we are done.
$\square$
\enddemo

Now let $\Cal X_1$ be the Kuranishi family
and let $\Cal X_2$ be obtained by pull-back by the transposition map. So the two families are locally equivalent by definition.
Now, by Lemma 6, they are not locally isomorphic.
\medskip
Observe that this trick gives only type (II) counterexamples.

\vskip .5cm
\subhead
2. Holomorphic families
\endsubhead

In this section, we prove

\proclaim{Theorem 3}
Let $\pi : \Cal X\to U$ be a holomorphic family of deformations. If $U$ is reduced and if
$h^0$ is constant in a neighborhood of $0$, then it has
the local isomorphism property.
\endproclaim

Notice the immediate Corollaries.

\proclaim{Corollary 3}
Let $X$ be a compact complex manifold such that $h^0$ is constant on its Kuranishi space $X$. Then any holomorphic deformation family 
of $X$ with reduced base has the local isomorphism property at $0$.
\endproclaim
 
\proclaim{Corollary 4}
Let $\pi : \Cal X\to U$ be a holomorphic family. Assume that $h^0$ is constant and that $U$ is reduced. Then every automorphism of $U$ lifts to an automorphism of $\Cal X$.
\endproclaim

\demo{Proof of Corollary 4}
Let $f$ be an automorphism of $U$, then $\Cal X$ and $f^*\Cal X$ are pointwise isomorphic.
So are locally isomorphic by Theorem 3. And this means that $f$ lifts. 
$\square$
\enddemo

On the other hand, recall that we gave in Lemma 6 an example of an automorphism of a reduced Kuranishi space which does not lift.

\demo{Proof of Theorem 3}
Assume $h^0$ constant in a neighborhood of $0$. Assume first that $\pi$ is a $1$-dimensional family parametrized
by the unit disk. Let $\pi' : \Cal X'\to\Bbb D$ be a pointwise isomorphic family. Let 
$$
f, g \ :\ \Bbb D\longrightarrow K
$$
such that $\Cal X=f^*\Cal K$ and $\Cal X'=g^*\Cal K$. We may assume without loss of generality that
\medskip

\noindent (i) The maps $f$ and $g$ are defined on the whole disk (otherwise shrink and uniformize).

\noindent (ii) The families $\Cal X$ and $\Cal X'$ are {\it equal} to $f^*\Cal K$ and $g^*\Cal K$, 
not only isomorphic (otherwise replace).
\medskip
We assume also without loss of generality that $K$ is reduced, since $f$ and $g$ map in fact onto $K^{\text{red}}$.
\medskip

Call $D$ the image of $f$, and $D'$ that of $g$. If $D$ or $D'$ is reduced to a point, then all the fibers of $\Cal X$ and of $\Cal X'$
are biholomorphic and Fischer-Grauert Theorem 
implies that {\it both} $D$ and $D'$ are reduced to a point. Both families are locally trivial,
hence locally isomorphic.
So we may assume
that $D$ and $D'$ are disks.

\medskip
Choose $(\phi_t)_{t\in\Bbb D}$ a family of pointwise biholomorphisms
$$
\phi_t \ :\ K_{f(t)} \longrightarrow K_{g(t)}
$$

\proclaim{Lemma 7}
There exists a dense subset of $\Bbb D$ such that, for each $t$ in this subset, there exists a sequence $(t_n)_{n\in\Bbb N}$ with all $t_n$ different from $t$ and with
$(\phi_{t_n})$ converging to $\phi_t$ in Sobolev norms as $n$ goes to infinity.
\endproclaim

\demo{Proof}
It is inspired from \cite{F-G}. Since the set $\Bbb D$ is uncountable,
the sequence $(\phi_t)_{t\in\Bbb D}$ contains an accumulation point for the Sobolev topology. This comes from the fact that $Diff(X^{\text{diff}})$ endowed with the Sobolev topology contains a countable dense sequence (compare with 
[Bourbaki, Topologie G\'en\'erale, Chapitre 10, Th\'eor\`eme 3.1] which gives a proof for the topology of uniform convergence). 
Moreover, given any open set $U'$ of $\Bbb D$ the same is true
for the subset $(\phi_t)_{t\in U'}$. Hence the claim.
$\square$
\enddemo

As a consequence, fix a neighborhood $V_0$ of $0$ in $K$. Then Lemma 7 implies that there exists a sequence 
$(t_n)_{n\in\Bbb N}\in V_0$ converging to some point $t_{\infty}$ of $V_0$ with $\phi_{t_n}$ converging to $\phi_{t_\infty}$ in
Sobolev norms. 
\medskip

Assume that $f(t_{\infty})$ is equal
to $g(t_{\infty})$ and that $\phi_{t_{\infty}}$ is the identity.   
By Proposition K2, assuming $V_0$ small enough, this would mean that we must have
$$
f(t_n)=g(t_n)
\leqno n\geq n_0
$$
for $n_0$ big enough. Now, this implies that $f-g$ is a holomorphic function on the disk with a non-discrete set of zeros, hence
$f\equiv g$ and we are done. Observe that $K$ is naturally embedded in the vector space $H^1(X,\Theta)$ as an analytic set,
hence the difference $f-g$ is meaningful as holomorphic map from $\Bbb D$ to $H^1(X,\Theta)$.
\medskip
In the general case, things become more complicated, but the previous pattern can be used as a guideline to proceed. Consider
the embedded family $\Cal K_{\vert D'}\to D'$ and write
$$
\Cal K_{\vert D'}=(D'\times X^{\text{diff}}, \alpha)
$$
where the complex operator $\alpha_t\in A^1$ turns $\{t\}\times X^{\text{diff}}$ into the complex manifold $K(t)$.
\medskip
Remark that we have a diffeomorphism 
$$
\phi_{t_{\infty}}^{-1} \ :\ K_{g(t_{\infty})}\longrightarrow  K_{f(t_{\infty})}
$$

Remark also that, by Corollary 2, the set $K$ is versal at both $f(t_{\infty})$ and $g(t_{\infty})$. So by Lemma 1, there exists
a diffeomorphism $(\Psi, \psi)$ of $\Cal K$ defined on a neighborhood of $f(t_{\infty})$ which extends $\phi_{t_{\infty}}^{-1}$. 
To simplify the notations, we identify in this proof a point $t$ of $K$ and the integrable almost-complex operator
$\alpha_t=G^{-1}(t)$ defining $K_t$. With this convention, $\psi$ is
constructed as a composition of
$$
(\phi_{t_{\infty}}^{-1})_*\ :\ A^1\longrightarrow A^1
$$
with the map $\Xi$ of Proposition K1. This gives us a new realization
$$
h\equiv \psi\circ g \ :\ U'\subset \Bbb D\longrightarrow K
$$
defined on a neighborhood $U'$ of $t_{\infty}$ such that $\Cal X'$ is locally isomorphic to $h^*K$ at $t_{\infty}$.
\medskip
But now we can make use of Proposition K2. The sequence
$$
\Psi\circ \phi_{t_n}\ :\ K_{f(t_n)}\longrightarrow K_{h(t_n)}
$$
converges in Sobolev norms to
$$
\Psi\circ \phi_{t_{\infty}}\ :\ K_{f(t_{\infty})}\aro >\text{Identity}>> K_{h(t_{\infty})}=K_{f(t_{\infty})}
$$
hence $f$ and $h$ take the same values not only at $t_\infty$ but also at every $t_n$ for $n$ big enough. 
Moreover, still by Proposition K2, and since we assumed that $K$ is reduced, the map $h$ is the unique map such that $\Cal X'$ is locally isomorphic to $h^*\Cal K$
at $t_{\infty}$ (provided $K$ is based and marked at $f(t_{\infty})$ and provided a marking of $\Cal X'$ is fixed at $t_{\infty}$ and
asked to be preserved). Since the family $\Cal X'$ is a holomorphic family, $h$ must be holomorphic. So as before we have
$f\equiv h$.

\remark{Remark}
This is just another way of saying that $K$ is universal with respect to families with $h^0$ constant equal to $h^0(0)$. Proposition K2
was proved by Kuranishi to have this type of result.
\endremark
\medskip

We claim that $\Cal X'$ is isomorphic to $h^*\Cal K$ over the whole disk $\Bbb D$, and not only over a neighborhood
of $t_{\infty}$ in $\Bbb D$. 
\medskip
This can be proven as follows. Let $U'\subset\Bbb D$ be the maximal subset of $\Bbb D$ such that
$\Cal X'_{\vert U'}$ is isomorphic to $h_{\vert U'}^*\Cal K$. Let $t\in\Bbb D$ and let $c$ be a path in $\Bbb D$ joining $c(0)=t_\infty$ to 
$c(1)=t$. We will prove that $t$ is in $U'$.
\medskip
The problem that could appear is that $(\phi_{t_{\infty}}^{-1})_*f(c)$, which is a path in $A^1$, 
is not fully included in the domain of definition $A$ of $\Xi$.
Let $K'\subset A^1$ be the Kuranishi space of $X'_{f(c(1))}$ based at $(\phi_{t_{\infty}}^{-1})_*f(c)$ (which is reduced since $K$ is reduced). Let $\Xi'$ be the map of Proposition
K1 defined in a neighborhood $A'$ of $(\phi_{t_{\infty}}^{-1})_*f(c(1))$ in $A^1$. For simplicity, assume that the whole path
$(\phi_{t_{\infty}}^{-1})_*f(c)$ is included in $A\cup A'$. Take a point $s\in [0,1]$ such that $(\phi_{t_{\infty}}^{-1})_*f(c(s))$
lies in the intersection of $A$ and $A'$. Then there exists a local isomorphism between the pointed analytic sets
$(K, \Xi((\phi_{t_{\infty}}^{-1})_*f(c(s)))$ and $(K',\Xi'((\phi_{t_{\infty}}^{-1})_*f(c(s)))$ since these two spaces
are versal for $X'_{f(c(s))}$. And this isomorphism can be chosen in such a way that the image of 
$\Xi((\phi_{t_{\infty}}^{-1})_*f(c))$ is sent to $\Xi'((\phi_{t_{\infty}}^{-1})_*f(c))$ in a neighborhood of $s$, still by the universality
property. 
\medskip

Let us sum up. We can glue $K$ and $K'$ to obtain an analytic space $\tilde K$ such that $h$ extends as $\tilde h$ 
along $c$ in such a way
that $\Cal X'$ is isomorphic to $\tilde h^*\tilde\Cal K$ along the full path $c$. Still by universality, in our case, $\tilde h$
must be equal to $h$, so that $t$ is in $U$. In particular, observe that the image of $\tilde h$ stays in $K\subset \tilde K$.
So the claim is proved.
\medskip
Now, we have
$$
\Cal X'\simeq h^*\Cal K=f^*\Cal K\simeq \Cal X
$$
on the whole disk (the symbol $\simeq$ meaning isomorphic). In other words, $\Cal X$ and $\Cal X'$ are locally isomorphic in a 
neighborhood of $0$. This 
proves the Theorem for $1$-dimensional families.
\medskip
Let us now assume that the families $\Cal X$ and $\Cal X'$ are $p$-dimensional.
We will use general arguments (already used in \cite{We}, though not exactly in the same way) 
to pass from the one-dimensional to the general case.
\medskip
By a Theorem of Namba \cite{Na, Theorem 2}, the union $\Cal H$ of pointwise holomorphic maps from $X_t$ to $X'_t$ for all $t$ can be
endowed with a structure of a reduced analytic space such that the natural projection map $p : \Cal H\to U$ is
holomorphic and surjective. Moreover, the topology of $\Cal H$ is that of uniform convergence.
\medskip
Let $\Cal S\subset\Cal H$ be the subset of pointwise isomorphisms. It is an open set of $\Cal H$ so a reduced analytic space with
a holomorphic (still surjective in this particular case) projection map $p$. This openness property can be shown as follows.
Given $\phi$, an isomorphism between $X_t$ and $X'_t$ for a fixed $t$, every $\psi$ close enough from $\phi$ in the topology of uniform
convergence is a local isomorphism at each point. We just have to prove now that $\psi$ must be bijective. Forgetting the complex
structures we can see $\phi$ and $\psi$ as maps of $X^{\text{diff}}$, using differentiable trivializations. Since $X^{\text{diff}}$ is compact and
$\psi$ locally bijective, $\psi$ is surjective. Besides,
still by compacity, there exists a finite open covering of $X^{\text{diff}}$ such that any map close enough from $\phi$ is injective when 
restricted to any member of this covering. Assume $\psi$ is not globally injective. Then, we could find a sequence of non-injective
maps $\psi_n$ converging 
uniformly onto $\phi$. So there would be two sequences of points $(x_n)$ and $(y_n)$ such that
$$
x_n\not =y_n\qquad \psi_n(x_n)=\psi_n(y_n)
\leqno n\in\Bbb N
$$

By compacity of $X^{\text{diff}}$, they will converge to some points $x$ and $y$ such that $\phi(x)=\phi(y)$, hence $x=y$. This clearly
contradicts the previous property of local injectivity of all $\phi_n$ on a fixed covering.
\medskip
To finish the proof of Theorem 3, it is enough to show that $p : S\to U$ has a local holomorphic section at $0$. But now, we conclude 
from what we did for $1$-dimensional families that $p$ has local holomorphic sections at $0$ along every embedded disk $\Bbb D$ in
$U$. Fix one of these local sections, say $\sigma$. Take another such section $\sigma'$. Then $\sigma$ and $\sigma'$ differ by
composition (at the source) by an automorphism of $X_0$ and by composition (at the target) by an automorphism of $X'_0$. If these automorphisms belong to the connected component of the identity, since
$U$ is reduced and $h^0$ is constant, both extend locally as automorphisms of the nearby fibers \cite{Gr}.
But this means exactly that, composing $\sigma'$ with these extensions, one may assume without loss of generality that
$\sigma'$ takes the same value at $0$ as $\sigma$. Using this trick and taking account that the number of connected components of the automorphism group of $X_0$ is countable, 
we see that there exist local sections with the same value at $0$
for almost every disk embedded in $U$ passing through $0$. 
\medskip
Now, by a Proposition of Grauert and Kerner \cite{G-K}, there exists an analytic embedding of a neighborhood $S$ of $\sigma(0)$ in 
$\Cal S$
$$
i\ :\ S\longrightarrow \Bbb D^{\dim p^{-1}(\{0\})}\times U
$$
such that the following diagram commutes
$$
\CD
S \aro > i>>\Bbb D^{\dim p^{-1}(\{0\})}\times U\\
\aro V p VV \aro VV\text{2nd proj.}V\\
U\aro >> \text{Identity}> U
\endCD
$$

Observe that the dimension of $p^{-1}(\{0\})$ is $h^0(0)$ and that, since $h^0$ is constant, $i(p^{-1}(s))$ is an open set
of $\Bbb D^{h^0(0)}\times\{s\}$. On the other hand, by what preceeds, $p(S)$ must be equal to an open neighborhood of $0$ in $U$
(because $U$ is reduced). As a consequence, $i$ is a local isomorphism, which exactly means that $\Cal X$ and $\Cal X'$ are locally
isomorphic at $0$.
$\square$
 \enddemo

\remark{Remark}
The last strategy (using Namba's Theorem and so on) cannot be used directly to obtain the result for $1$-dimensional families.
Indeed, it is not possible to exclude the case of $p^{-1}(\{0\})$ being isolated from the other fibers, so that in the diagram above,
the image $p(S)$ reduces to $0$. The only fact that can be proven directly is that, if we know that there exists a sequence
$$
\phi_{t_n}\ :\ X_{t_n}\longrightarrow X'_{t_n}
$$
converging uniformly to some $\phi_0: X_0\to X'_0$, then the two families are locally isomorphic at $0$. This is just because,
in this case, $p(S)$ is an analytic set of $\Bbb D$ (we are in the $1$-dimensional case) containing 
an infinite sequence $(t_n)$ accumulating on $0$. So $p(S)$ must contain an open neighborhood of $0$. Now, we obtained the same 
conclusion using Proposition K2.
\endremark

\remark{Remark}
In the non-reduced case, Theorem 3 is false, as shown by the following easy example. Consider the upper half-plane $\Bbb H$ of $\Bbb C$
as the parameter space of elliptic curves. Let $\Cal H\to \Bbb H$ the versal (at each point) associated family.
Choose a point $\tau\in\Bbb H$. Take $U$ to be the double point 
$$
U=\{t^2=0\quad\vert\quad t\in\Bbb C\}
$$ 
Let $\pi : \Cal X_1\to U$ be the constant family obtained by pull-back by a constant map from $U$ to $\Bbb H$ (with value $\tau$).
Now, since $U$ is not reduced, there exists also non-constant morphisms from $U$ to $\Bbb H$. Let $f$ be the unique such morphism
sending the single point of $U$ to $\tau$ and the vector $\partial/\partial t$ of the Zariski tangent space of $U$ to the horizontal unit
vector of $\Bbb H$ based at $\tau$. Define $\Cal X_2$ as $f^*\Bbb H$. Then $\Cal X_1$ and $\Cal X_2$ are obviously pointwise isomorphic,
but they are not locally isomorphic, by computation of their Kodaira-Spencer map. It is $0$ for $\Cal X_1$, and not zero for $\Cal X_2$.
\endremark

\vskip.5cm

\subhead
3. Type (II)-counterexamples to the equivalence problem
\endsubhead

We derive now a characterization of type (II) counterexamples to the equivalence problem in the one-dimensional case.

\proclaim{Theorem 4}
The following statements are equivalent.
\medskip
\noindent (i) The one-dimensional holomorphic families $\pi:\Cal X\to\Bbb D$ and $\pi':\Cal X'\to\Bbb D$ form a type (II) counterexample
to the equivalence problem.

\noindent (ii) Both are obtained from the same jumping family $\pi'' : \Cal X''\to\Bbb D$
by pull-backs by some maps. Moreover, the degrees of these maps (as ramified coverings of $\Bbb D$) are different.
\endproclaim

\demo{Proof}
Assume that $\Cal X$ (respectively $\Cal X'$) are obtained from the Kuranishi space of $X$ by pull-back by some map $f$ 
(respectively $h$). Call $D$ the image of $f$ and $D'$ that of $h$.
Shrinking the domains of definition if necessary to have the same image $D''\subset D\cap D'$ 
and uniformizing at the source and at the target by unit disks,
we obtain the following diagram
$$
\CD
\Bbb D &\aro >\text{uniform.} >>f^{-1}(D'')\subset \Bbb D& \aro > f>> D''&\aro >\text{uniform.}>> \Bbb D&\\
&&&&&&\aro V \text{Id} VV &\aro VV \text{Id} V\\
\Bbb D &\aro >\text{uniform.} >> h^{-1}(D'')\subset\Bbb D&\aro >h >> D''&\aro >\text{uniform.}>>\Bbb D&
\endCD
$$

To simplify, we still denote by $f$ (respectively by $h$) the composition of the top arrows (respectively of the bottom arrows).
Moreover, we denote by $\pi'' : \Cal X''\to\Bbb D$ the target family
and replace $\Cal X$ (respectively $\Cal X'$) by $f^*\Cal X''$ (respectively $h^*\Cal X''$).

\medskip
We may assume without loss of generality that $f$ and $g$ are unramified coverings over $\Bbb D^*$ of respective 
degrees $n$ and $m$. So we have \cite{Fo, Theorem 5.11}
$$
f(z)=z^n\qquad g(z)=z^m
\leqno z\in\Bbb D
$$
changing the uniformizing maps at the source by a rotation if necessary.
\medskip

If $m$ and $n$ are equal, then $\Cal X$ and $\Cal X'$ are locally equivalent at $0$. So 
assume $n>m$.
\medskip

Now, from the one hand, by definition of the pull-back, for all $t\in\Bbb D$, 
the fibers $X_t$ and $X''_{t^n}$ are 
biholomorphic, as well as $X'_t$ and $X''_{t^m}$. 
And from the other hand, the assumption for the families of being pointwise isomorphic means in this new setting that
there exists
$$
\Phi\ :\ (U_1\subset \Bbb D, 0)\longrightarrow (U_2\subset \Bbb D, 0)
$$
a biholomorphism such that $X_t$ and $X'_{\phi(t)}$ are biholomorphic. Hence by transitivity,
$X_t$ and $X_{\phi^{-1}(t^{n/m})}$ are biholomorphic for every choice of a determination of $t^{n/m}$. Observe that this is valid
for $t$ belonging to a sufficiently small neighborhood $U'_1$ of $0$ in $\Bbb D$.
Set
$$
C=\{t\in\Bbb D\quad\vert\quad \vert t\vert=\lambda\}
$$
for $\lambda$ a fixed real number, which is supposed small enough to have $C\subset U'_1$. 

\proclaim{Lemma 8}
Let $t_0\in C$. Then the closure of the set
$$
E_{t_0}=\{t\in \Bbb D\quad\vert\quad X_{t_0}\text{ is biholomorphic to } X_t\}
$$
contains $C$.
\endproclaim

\demo{Proof of Lemma 8}
Assume first that $\phi$ is equal to $a\cdot \text{Id}$ for $a$ non-zero. Defining $\alpha_k$ for $k\in\Bbb N$ by induction 
$$
\left\{
\eqalign{
\alpha_0=&a^{-1}\cr
\alpha_{k+1}&=\alpha^{-1}\cdot \alpha_k^{n/m}
}
\right .
$$
(we choose a determination of $\alpha_k\mapsto \alpha_k^{n/m}$ for each $k$), we have that $X_t$ and $X_{\alpha_k\cdot (t^{n/m})^k}$ are biholomorphic.
In particular, all the points of
$$
\{t_0\exp (2i\pi l(m/n)^k)\quad\vert\quad k>0,\ l\in\Bbb Z\}
$$
correspond to $X_{t_0}$ proving the density of $E_{t_0}$ in $C$.
\medskip
Now, if $\phi$ is not a homothety, it admits a Taylor expansion
$$
\phi(t)=at+\text{ higher order terms}
$$

Besides, one has that $t\in E_{t_0}$ as soon as $t^{n/m}=t_0^{n/m}$, or
$$
(\phi^{-1}(t^{n/m}))^{n/m}=(\phi^{-1}(t_0^{n/m}))^{n/m}
$$
or more generally
$$
\left (\phi^{-1}\bigl (\hdots (\phi^{-1}(t^{n/m}))^{n/m}\hdots \bigl )^{n/m}\right )^{n/m}=
\left (\phi^{-1}\bigl (\hdots (\phi^{-1}(t_0^{n/m}))^{n/m}\hdots \bigl )^{n/m}\right )^{n/m}
$$

Using the Taylor expansion of $\phi$ together with the fact that $n/m>1$, we obtain that the sequence
$$
\dfrac{\left (\phi^{-1}\bigl (\hdots (\phi^{-1}(t^{n/m}))^{n/m}\hdots \bigl )^{n/m}\right )^{n/m}}
{\alpha_k\cdot (t^{n/m})^k}
$$
tends to $1$ as $k$ goes to infinity. In this expression, the determinations of the $n/m$-th power are chosen at each step according to
the choices made for $\alpha_k$.
\medskip
This means that, given
$$
t=t_0\exp (2i\pi l(m/n)^k)
$$
for some fixed $k>0$ and $l\in\Bbb Z$, and given any $\epsilon >0$, there exists $t'\in \Bbb D$ which is $\epsilon$-close to $t$ such that
$t'$ belongs to $E_{t_0}$. This is enough to conclude that the closure of $E_{t_0}$ contains $C$.
$\square$
\enddemo

Hence, there exists a dense subset of points corresponding to $X_{t_0}$ in any annulus around the circle $\vert z\vert =\vert t_0\vert$. 
Following \cite{Gr}, the function
$h^1$ is constant on a Zariski open subset of $\Bbb D$. So we may assume that it is constant on $\Bbb D^*$. That means that the 
differentiable family of deformations parametrized by $\vert z\vert =\vert t_0\vert$ is a regular one. 

\proclaim{Proposition 3}
All points of the circle $C=\{\vert z\vert=\vert t_0\vert\}$ in $\Bbb D$ correspond to the same complex manifold $X_{t_0}$.
\endproclaim

\demo{Proof}
This is a  step by step adaptation of  
the proof of \cite{F-G}. We will prove that the Kodaira-Spencer map is zero for a dense subset of points, hence by regularity for all
points, so that all points correspond to the same manifold $X_{t_0}$.
We will explain in details how to modify the proof of \cite{F-G} so that it generalizes to this case, but will refer freely to \cite{F-G} 
for the common parts.
\medskip
Choose $\epsilon>0$. Choose also an annulus $A$ around $C$. Choose finally  a differentiable
trivialization 
$$
T\ :\ \pi^{-1}(\{s\in A\quad\vert\quad \vert s-t_0\vert <\epsilon\})\longrightarrow X_{t_0}
$$
with $T_{t_0}\equiv \text{Id}$. 
\medskip
For every $t$ in $C$, define a diffeomorphism $\tilde\alpha_t$ from $X_{t_0}$ to $X_t$ as follows. First, choose some $t'\in A$ such that
\bigskip
\noindent (i) We have $\vert t'-t_0\vert <\min (\vert t-t_0\vert, \epsilon)$.

\noindent (ii) The parameter $t'$ belongs to the set 
$$
E_{t}\cap \pi^{-1}(\{s\in A\quad\vert\quad \vert s-t_0\vert <\epsilon\})
$$
\bigskip

This is possible by Lemma 8. By what preceeds, there exists a biholomorphism $\beta_t$ between $X_{t'}$ and $X_t$.
Define 
$$
\tilde\alpha_t\equiv \beta_t\circ T^{-1}_{t'}
$$

First notice that the set
$$
E=\{t\in C\quad\vert\quad \exists (t_n)_{n\in\Bbb N}\in C\text{ such that } (\tilde\alpha_{t_n}) 
\text{ uniformly converges to } \tilde\alpha_t\}
$$
is dense in $C$. This comes from the fact that the set of continuous maps from $X$ to $\pi^{-1}(C)$ is of countable type, see
\cite{F-G} and the appendix.
\medskip
Let $t\in E$. Without loss of generality, we may assume that there exists a finite set of submersion charts 
$$
\CD
U_i\subset \pi^{-1}(\{s\in A\quad\vert\quad \vert s-t\vert <\epsilon\}) \aro >\psi_i>> \Bbb C^{\dim X}\times 
\{s\in A\quad\vert\quad \vert s-t\vert <\epsilon\} \\
\aro V\pi VV \aro VV \text{2nd projection}V\\
\{s\in A\quad\vert\quad \vert s-t\vert <\epsilon\}\aro >>\text{Identity} >\{s\in A\quad\vert\quad \vert s-t\vert <\epsilon\}
\endCD
$$
covering $\pi^{-1}(\{s\in A\quad\vert\quad \vert s-t\vert <\epsilon\})$.
\medskip
Set $\alpha_n\equiv \tilde\alpha_{t_n}\circ \alpha_t^{-1}$. Then the sequence $(\alpha_n)$ converges uniformly to the identity
of $X_t$. Let $(V_j)$ be a covering of $\pi^{-1}(\{s\in A\quad\vert\quad \vert s-t\vert <\epsilon\})$ by relatively compact
open sets with smooth boundaries such that there exists a refining map $r$ and an integer $n_0$ satisfying
$$
\forall n\geq n_0,\quad \alpha_n(V_j)\subset U_{r(j)}
$$
First, assume for simplicity that $h^0(t)=0$.
Let $x\in V_i\cap X_t$ and let $(z,t)$ be the coordinates in the chart $\psi_{r(i)}$. For $n\geq n_0$, define 
$$
\eqalign{
\xi_i^n(x)&=(\psi_{r(i)}^{-1})_*(\psi_{r(i)}\circ\alpha_n(x)-\psi_{r(i)}(x))\cr
&=(\psi_{r(i)}^{-1})_*(z(\alpha_n(x))-z(x), t(\alpha_n(x))-t(x))\cr
&=(\psi_{r(i)}^{-1})_*(z(\alpha_n(x))-z(x), t_n-t)
}
$$
where $(\psi_{r(i)}^{-1})_*$ denotes the pushforward of a vector field by the differential of $\psi_{r(i)}^{-1}$. 
This gives a smooth vector field on $V_i\cap X_t$ which is transverse to $X_t$ (since the $t$-coordinate is non-zero). Now, let
$$
M_n=\max_i\sup_{x\in V_i\cap X_t} \Vert (\psi_{r(i)})_*\xi_i^n(x)\Vert
$$ 
for some choice of a norm on $\Bbb C^{\dim X_t}\times\Bbb R$. Notice that $M_n$ is positive since it is bigger than $\vert t_n-t\vert$;
and that it is finite because of the finiteness of the number of charts and because of the relative compactness of the $V_i$.

\proclaim{Lemma 9}
The sequence $1/M_n(\xi_i^n)$ converges uniformly to a holomorphic vector field $\xi_i$ on $V_i\cap X_t$.
\endproclaim

\demo{Proof of Lemma 9}
Let $\eta_i^n=1/M_n(\psi_{r(i)})_*(\xi_i^n)$. It is a uniformly bounded sequence of functions on $D_i=\psi_{r(i)}(V_i\cap X_t)$. If we prove
that is an equicontinuous sequence, then Ascoli's Theorem will ensure the uniform convergence.
\medskip
Now, for all $j$ between $1$ and $\dim X$, the sequence 
$$
\bar\partial_j \eta_i^n\equiv \dfrac{\partial}{\partial \bar z_j}\eta_i^n
$$ 
is uniformly convergent to zero since we have
$$
\bar\partial_j \eta_i^n=\dfrac{1}{M_n} (\bar\partial_j (z(\alpha_n\circ\psi_{r(i)}^{-1}))-\bar\partial _j z, 
\bar\partial_j (t(\alpha_n\circ\psi_{r(i)}^{-1})))
$$
and since $\alpha_n$ tends uniformly to the identity, hence $z(\alpha_n\circ\psi_{r(i)}^{-1})$ tends uniformly to $z$
and $t(\alpha_n\circ\psi_{r(i)}^{-1}))$ to $t$.
\medskip
On the other hand, deriving with respect to $z_j$ the Bochner-Martinelli formula for $\eta_i^n$, one obtains, for $k=\dim X_t$,

$$
\displaylines{
\partial _j\eta_i^n(z)=\dfrac{k!}{(2i\pi)^k} \Bigl (\int_{\partial D_i} \sum_{\nu=1}^k ((-1)^{\nu }\eta_i^n)(\zeta)
\dfrac{(\bar \zeta_\nu-\bar z_\nu)(\bar\zeta_j-\bar z_j)}{\vert \zeta-z\vert ^{3k+1}}d\bar\zeta [\nu]\wedge d\zeta \hfill\cr
\hfill +\int_{D_i} \sum_{\nu=1}^k ((-1)^{\nu }\bar\partial_{\nu} \eta_i^n)(\zeta)
\dfrac{(\bar \zeta_\nu-\bar z_\nu)(\bar\zeta_j-\bar z_j)}{\vert \zeta-z\vert ^{3k+1}}d\bar\zeta [\nu]\wedge d\zeta \Bigl )
}
$$
where $d\bar\zeta[\nu]=d\bar\zeta_1\wedge\hdots\wedge d\bar\zeta_{\nu -1}\wedge d\bar\zeta_{\nu +1}\wedge \hdots\wedge d\bar\zeta_n$.

Since $(\eta_i^n)$ is a uniformly bounded sequence and since $(\bar\partial _j\eta_i^n)$ is uniformly convergent to zero, we obtain that 
$(\partial_j \eta_i^n)$ is also a uniformly bounded sequence. So $(\eta_i^n)$ is Lipschitz with a Lipschitz constant independant of $n$.
Therefore it is equicontinuous.
\medskip
Finally, since $(\bar\partial _j\eta_i^n)$ is uniformly convergent to zero, the limit is automatically holomorphic.
$\square$
\enddemo

Following \cite{F-G}, it is easy to prove that these $\xi_i$ glue together to define a global non-zero holomorphic vector field
$\xi$ on $X_t$. This vector field must be transverse to $X_t$ for we assumed $h^0(t)=0$. Hence, the Kodaira-Spencer map at $t$ is
zero.
\medskip
If $h^0(t)$ is not zero, one has first 
to modify each $\alpha_n$ by composition with a finite number of well-chosen automorphisms of $X_t$.
The construction of the holomorphic vector field $\xi$ is then exactly the same. And finally one uses the special properties of
this new sequence of $(\alpha_n)$ to prove that $\xi$ cannot be tangent to $X_t$. Details are exactly the same as in \cite{F-G}.
\medskip
As a consequence, one obtains that the family over $C$ has zero Kodaira-Spencer map on a dense subset of points, hence on $C$ as it is 
a regular family. And applying Theorem 6.2 of \cite{K-S1}, one has that this family is locally trivial at every point. 
In particular, all the
fibers correspond to the same compact complex manifold up to biholomorphism.
$\square$
\enddemo

But the existence of such a circle of biholomorphic fibers forces the foliation of Section III to be non-trivial.
From the previous proof, we deduce that all the circles $z=\vert t\vert$ of $\Cal X$
correspond to a unique complex structure, say $X_{t}$. 
Fix such a $t$ different from $0$.
This implies that the intersection of $\Bbb D$ with 
the leaf of the foliation passing through $t$ contains a circle of points. 
Since the foliation is holomorphic, this means that a neighborhood of this circle
corresponds to $X_t$. Let $s$ be in the boundary of this neighborhood. Then the same argument shows that a neighborhood
of the circle $\vert z\vert =\vert s\vert$ lies in the leaf passing through $s$. Now, the two previous neighborhoods must have non-empty
intersection which implies that $X_s$ and $X_t$ are biholomorphic.
\medskip
We conclude from that that all the points of $\Bbb D^*$ correspond to $X_t$. Hence, by Fischer-Grauert Theorem, $\Cal X''$ must
be a jumping family. 
\medskip
To prove the converse, we need to refine the argument given in the proof of Proposition 2. Consider the {\it local} Kodaira-Spencer
map of $\Cal X$ at $0$
$$
H^0(U,\Theta)\aro >\rho_{\Cal X} >> H^1(\Cal X_{\vert U},\Theta)
\leqno 0\in U\subset \Bbb D
$$
which represents the obstruction to lifting a holomorphic vector field in the base $U\subset \Bbb D$ to the family 
$\Cal X_{\vert U}=\pi^{-1}(U)$. The direct limit of $\rho_{\Cal X}$ for $U$ smaller and smaller 
gives the pointwise Kodaira-Spencer map used in the
proof of Proposition 2 and which represents the pointwise first-order obstruction to this lifting problem.
\medskip
But we can also define a pointwise $(p+1)$-th order obstruction for any $p\in\Bbb N$ and any $\xi\in H^0(U,\Theta)$ 
by taking the $p$-jet of $\rho_{\Cal X}(\xi)$ at $0$ 
(jet as local sections of $\Theta$) and passing to the direct limit. This defines a $(p+1)$-th order Kodaira-Spencer map
$$
J^p_0(T\Bbb D) \aro > \rho_{\Cal X}^{(p)}>> H^1(X_0,\Theta^{(p)})
$$
where $J^p_0(T\Bbb D)$ is the vector space of $p$-jets at $0$ of holomorphic vector fields of $\Bbb D$ defined in a neighborhood
of $0$ and $\Theta^{(p)}$ is the bundle of $p$-jets of holomorphic sections of $\Theta$ (cf \cite{Wa}).
\medskip
Since the local Kodaira-Spencer map satisfies a chain-rule property, so does $\rho_{\Cal X}^{(p)}$, Hence, starting from $\Cal X$, 
$\Cal X'$ pull-backs of $\Cal X''$ by maps $f$ and $g$, we obtain the following equality
$$
\rho_{\Cal X}^{(p)}\left (\dfrac{\partial}{\partial t}\right )=\rho_{\Cal X''}^{(p)}\left (f_*(j^p_o(\dfrac{\partial}{\partial t}))\right )
\qquad\text{ and }\qquad
\rho_{\Cal X'}^{(p)}\left (\dfrac{\partial}{\partial t}\right )=\rho_{\Cal X''}^{(p)}\left (g_*(j^p_0((\dfrac{\partial}{\partial t}))\right )
$$
with $f_*$ (respectively $g_*$) denoting the action of $f$ (respectively $g$) on $p$-jets of vector fields.
Now if $f$ has degree $n$ and $g$ degree $m$, the above $(p+1)$-th obstruction of $\Cal X$ vanishes for $p<n$ and does not
vanish for $p=n$, whereas that of $\Cal X'$
vanishes for $p<m$ and does not vanish for $p=m$. 
Hence, if $m$ and $n$ are different, the families $\Cal X$ and $\Cal X'$ are not locally isomorphic at $0$.
$\square$
\enddemo

Of course, this is no more true for higher-dimensional families. Starting from two pointwise isomorphic but not locally isomorphic
one-dimensional jumping families, one can take their products with a fixed family and obtain type (II) counterexamples which are
not coming from jumping families.

\vskip .5cm
\subhead
4. Differentiable families
\endsubhead

Things are completely different for differentiable families. In fact, we have

\proclaim{Theorem 5}

\noindent (i) Let $\pi :\Cal X\to V$ be a real analytic family. If $h^0$ is constant along the family, then it has the local isomorphism 
property.

\noindent (ii)
Some differentiable families $\pi :\Cal X\to I$ of 2-dimensional compact complex tori do not have the local isomorphism property.

Moreover, there exist counterexamples of type (I) among families of 2-di\-men\-sio\-nal compact complex tori.
\endproclaim

\demo{Proof}

\noindent (i) This is exactly the same proof as that of Theorem 3. For the $1$-dimensional part, 
we observe that the only properties of holomorphic maps used are properties of analytic functions. For the passage to higher
dimension, it is enough to embed pointwise isomorphic families $\Cal X$ and $\Cal X'$ in holomorphic
families $\Cal X_{\Bbb C}\to U$ and $\Cal X'_{\Bbb C}\to U$ with constant $h^0$. 
For example, one may take for $U$ the reduction
of the stratum $K^{\text{max}}$. Then the only difference is that the map $p : \Cal S\to U$ given by Namba's Theorem may not be surjective.
But the same argument shows that it has a holomorphic section at $0$ defined on an analytic subspace of $U$ containing $V$.

\remark{Remark}
The same proof shows that if two differentiable families over $V$ are pointwise isomorphic {\it and} locally isomorphic along
each path of $V$ containing $0$, then they are locally isomorphic.
\endremark
\medskip

\noindent (ii) Because of (i), a smooth family of tori not having the local isomorphism property at $0$ must be flat
at $0$.
\medskip

Recall \cite{K-S2} that the open set
$$
M=\{A\in M_2(\Bbb C)\quad\vert\quad \det (\Im A)>0\}
$$
is a versal (and even universal) deformation space for every 2-dimensional compact complex torus. A point $A=(A_1,A_2)$ of $M$
corresponds to the quotient of $\Bbb C^2$ by the lattice generated by
$$
(1,0)\quad (0,1)\quad A_1\quad A_2
$$
Notice that every torus can be obtained as such a quotient. Two different points $A$ and $B$
of $M$ define the same torus up to biholomorphism if
and only if there exists
$$
\gamma=\pmatrix
\gamma_{11} &\gamma_{12}\\
\gamma_{21} &\gamma_{22}
\endpmatrix
\in SL_4(\Bbb Z)\quad\text{ such that }\quad
B=A\cdot\gamma=(\gamma_{11}+A\gamma_{21})^{-1}(\gamma_{12}+A\gamma_{22})
$$
Finally, $h^0$ is constant equal to 4 (given by the translations), so the condition of Theorem 3 is satisfied.
\medskip

Let 
$$
\Omega_0=\pmatrix 
i &0 \\
0 & i
\endpmatrix
\quad\text{ and }\quad
\Omega(t)=\pmatrix
i+t & b(t)\\
c(t) & i+t
\endpmatrix
\eqno t\in \Bbb R
$$
and let $X_t$ be the corresponding tori. The smooth functions $b$ and $c$ satisfy
\bigskip
\noindent 1. They are smoothly flat at zero, i.e. all their derivatives at zero are zero.

\noindent 2. We have $b(0)=c(0)=0$ and $b(t)>c(t)>0$ for $t$ different from zero.
\bigskip

The path $\Omega$ in $M$ defines a differentiable family of 2-dimensional compact complex tori centered at $X_0$.
Define $\Omega_1\equiv \Omega$ and
$$
\Omega_2(t)=\left\{
\eqalign{
&\Omega_1(t) \quad\text{ if } t\leq 0\cr
& ^t \Omega_1(t) \quad\text{ if } t\geq 0
}
\right .
$$

Remark that conditions 1 and 2 imply that $\Omega_2$ is also a smooth path.
\medskip

We claim that the corresponding families $\Cal X_1\to\Omega_1$ and $\Cal X_2\to\Omega_2$ are pointwise isomorphic but not locally 
isomorphic at $0$.
\medskip
First note that, for all $t$,
$$
^t\Omega_1(t)=\pmatrix 0 &1\\ 1 &0\endpmatrix \cdot\Omega_1(t)\cdot \pmatrix 0 &1\\ 1 &0\endpmatrix=\Omega_1(t)\cdot\gamma
$$
for
$$
\gamma=\pmatrix
0 &1 &0 &0\\
1 &0 &0 &0\\
0 &0 &0 &1\\
0 &0 &1 &0
\endpmatrix
\in SL_4(\Bbb Z)
$$
That implies that, for $t>0$, the map
$$
(z,w)\in\Bbb C^2\longmapsto (w,z)\in\Bbb C^2
$$
descends as a biholomorphism between $X_1(t)$ and $X_2(t)$. So the families are pointwise isomorphic.
\medskip
On the other hand, for a generic lattice, it is well-known that the automorphism group of a torus is generated by translations and by
$-\text{Id}$. Indeed, for this particular choice of matrices $\Omega(t)$, it is straightforward that this is the case if the numbers
$i+t$, $b(t)$, $c(t)$, their squares and all the products of two of them are linearly independent over $\Bbb Q$. 
Hence, for generic $t$, the tori $X_1(t)$ and $X_2(t)$ have no other automorphisms than these ones. This allows to find
sequences $(t'_n)_{n\in\Bbb N}$ of negative numbers and $(t''_n)_{n\in\Bbb N}$ of positive numbers converging to $0$ such that
\bigskip
\noindent (i) For each $n$, up to translations, the only biholomorphisms between
$X_1(t'_n)$ and $X_2(t'_n)$ are the projection of $\pm \text{Id}$ on $\Bbb C^2$.

\noindent (ii) For each $n$, up to translations, the only biholomorphisms between
$X_1(t''_n)$ and $X_2(t''_n)$ are the projection of $\pm \pmatrix 0 &1\\ 1 &0\endpmatrix$ on $\Bbb C^2$.

\bigskip

Suppose now that $\Cal X_1$ and $\Cal X_2$ are locally isomorphic at $0$. Then there would exist a family $(\Phi_t)$ of biholomorphisms
of $\Bbb C^2$ (for $t$ in a neighborhood of $0$) such that
\bigskip
\noindent (i) It is smooth in $t$.

\noindent (ii) Every $\Phi_t$ descends as a biholomorphism between $X_1(t)$ and $X_2(t)$.
\bigskip
 
But, by what preceeds, at $t'_n$ the map $\Phi_t$ must be $\pm \text{Id}$ up to a translation factor, whereas at $t''_n$, it must be
$\pm \pmatrix 0 &1\\ 1 &0\endpmatrix$ up to a translation factor. Since these two sequences do not converge to the same type of limit
when $n$ goes to infinity, we arrive to a contradiction. The families $\Omega_1$ and $\Omega_2$ are not
locally isomorphic at $0$.
\medskip

On the other hand, the previous family still has the local isomorphism property when restricted to $(-\infty, 0]$ and $[0,\infty)$.
Nevertheless, it is easy to modify it in order to have a counterexample even when restricted to $(-\infty, 0]$ and $[0,\infty)$.
\medskip
Start with the same path $\Omega$ as before, but this time assume that the functions $b$ and $c$ satisfy
\bigskip
\noindent 1. There exists a sequence $(t_n)_{n\in\Bbb N}$ of positive numbers converging to $0$ such that $b$ and $c$ are zero 
and flat at all $t_n$.

\noindent 2. We have $b$ and $c$ even.

\noindent 3. We have $b(t)\not = c(t)$ for $t$ positive and not belonging to the sequence $(t_n)$.
\bigskip

For example, let
$$
h(t)=\left \{\eqalign
{
&0 \quad\text{ for }t\leq 0 \cr
&\exp (-1/t)\quad\text{ otherwise }
}
\right .
$$
and
$$
f : t\in \Bbb R\longmapsto \sum_{p\in\Bbb Z}h(t+p)\cdot h(-t-p+1)\in\Bbb R
$$
and finally
$$
b\equiv \alpha h(\vert -\vert )\cdot f(\log \vert -\vert )\qquad b\equiv \beta h (\vert -\vert )\cdot f(\log \vert -\vert )
$$
for $\beta\not = \alpha$. In this case, we have $(t_n)=\exp (-n)$.
\medskip
The path $\Omega$ in $M$ defines a differentiable family of 2-dimensional compact complex tori centered at $X_0$.
Define $\Omega_1\equiv \Omega$ and
$$
\Omega_2(t)=\left\{
\eqalign{
&\Omega_1(t) \quad\text{ if } \vert t\vert\in [t_{2n},t_{2n+1}]\text{ for some }n\cr
& ^t \Omega_1(t) \quad\text{ if } \vert t\vert\in [t_{2n-1},t_{2n}]\text{ for some }n
}
\right .
$$
That implies that, for $t\in [t_{2n-1},t_{2n}]\text{ for some }n$, the map
$$
(z,w)\in\Bbb C^2\longmapsto (w,z)\in\Bbb C^2
$$
descends as a biholomorphism between $X_1(t)$ and $X_2(t)$. In particular, it defines an automorphism of $X_0$ and of $X_1(t_n)=X_2(t_n)$
for all $n$. This proves the pointwise isomorphism between the fibers.
\medskip
On the other hand, as in the previous example, one can find
sequences $(t'_n)_{n\in\Bbb N}$ and $(t''_n)_{n\in\Bbb N}$ of positive numbers converging to $0$ such that
\bigskip
\noindent (i) For each $n$, we have $t'_n\in [t_{2n},t_{2n+1}]$ and, up to translations, the only biholomorphisms between
$X_1(t'_n)$ and $X_2(t'_n)$ are the projection of $\pm \text{Id}$ on $\Bbb C^2$.

\noindent (ii) For each $n$, we have $t''_n\in [t_{2n-1},t_{2n}]$ and, up to translations, the only biholomorphisms between
$X_1(t''_n)$ and $X_2(t''_n)$ are the projection of $\pm \pmatrix 0 &1\\ 1 &0\endpmatrix$ on $\Bbb C^2$.
\bigskip

This is enough to prove that these two families, when seen as families over $[0,\infty)$, are not locally isomorphic at $0$.
Since the functions $b$ and $c$ are even, the same is true over $(-\infty, 0]$.
$\square$
\enddemo

In the differentiable case, it seems difficult to give a sufficient condition to have the local isomorphism property, except for the following trivial one.

\proclaim{Proposition 4}
Let $X$ be a compact complex manifold. Suppose that $K$ is a local moduli space for $X$
(that means that two different points of $X$ corresponds to two non-biholomorphic manifolds). 
Then every holomorphic (over a reduced base) 
as well as differentiable deformation family of $X$ has the local isomorphism property.
\endproclaim

\demo{Proof}
In this case, given any deformation family $\Cal X$ of $X$, the map from the parameter space of $\Cal X$ to $K$ is uniquely
determined by the pointwise complex structure of the fibers.
$\square$
\enddemo

\vskip .5cm
\subhead
5. Universality
\endsubhead

Let us finish this section by a comparison between our uniformization problems and the problem of universality of the Kuranishi
space.

\proclaim{Proposition 5}
Let $K$ be the Kuranishi space of some compact complex manifold $X$. Then the following statements are equivalent.
\medskip
\noindent (i) The space $K$ is universal for differentiable families.

\noindent (ii) The space $K$ is universal for holomorphic families over a reduced base.

\noindent (iii) The foliation of $K$ described in Section III is trivial.

\noindent (iv) The function $h^0$ is constant on $K$.
\endproclaim

\demo{Proof}
The equivalence $(iii)\iff (iv)$ is given by Theorem 2. 
The implications $(iv)\Rightarrow (i)$ and $(iv)\Rightarrow (ii)$ are immediate consequences of Proposition K2. Indeed, it is used in \cite{Ku2} to prove that.
The converse $(ii)\Rightarrow (iv)$ is proved in \cite{Wa} and \cite{Wa2}. Indeed, (ii) can be replaced by: the space $K^{\text{red}}$ 
is universal for holomorphic families over a reduced base. Now, from \cite{Wa, Proposition 4.2}, we have that it is the case if and only a certain extension problem
(called the second extension problem in \cite{Wa}) is solvable for $K^{\text{red}}$. Then one uses \cite{Wa2, p. 349} to conclude. 
\medskip
So we just need to prove $(i)\Rightarrow (iv)$. Assume $h^0$ non-constant. Then, by Proposition 1, there exists an automorphism $\phi$ of $X$ isotopic to
the identity such that any extension as a diffeomorphism of $\Cal K$ does not project onto the identity of $K$. Let $\Phi$
be such an extension. Still by Proposition 1, we may assume that $\Phi$ is isotopic to the identity. Let $(\Phi_t)_{t\in [0,1]}$
be such an isotopy, $\Phi_0$ being the identity map. 
Now set $\Psi(-)=\Phi_{\lambda (-)}(-)$, for some smooth function $\lambda : K\to [0,1]$ satisfying
\bigskip
\noindent (i) $\lambda_0(0)=0$.

\noindent (ii) $\det \text{Jac}_0\ \lambda\not =0$.
\bigskip
For a good choice of $\lambda$, the map $\Psi$ is a local diffeomorphism at $0$. Indeed, a direct computation shows that
$$
\text{Jac}_0\ \Psi=\text{Id}+\text{Jac}_0\ \lambda\cdot \dfrac{\partial \Phi_t}{\partial t}_{\vert t=0}
$$
so it is enough to take $\Vert \text{Jac}_0 \lambda\Vert$ very small.
\medskip
Recall now that $\Phi_t$ may be chosen so that, for all $t$, its germ at $0$ does not project as the germ of the identity (see the
proof of Lemma 2 and Proposition 1). From this, we deduce that the germ of $\Psi$ at $X_0$ is not the identity, even if
$\Psi_{\vert X_0}$ is the identity of $X_0$. In other words, one can
find a path $c$ in $K$ passing through $0$ whose image by $\Psi$ is different from $c$. But this means that the family corresponding
to $c$ is locally isomorphic to the family corresponding to $\Psi(c)$, with the same identification at $0$. Hence $K$ is not universal
for differentiable families.
$\square$
\enddemo

Consider now the case where $K$ is non-reduced. For example, assume

$$
K=\{t^2=0\quad\vert\quad t\in\Bbb C\}
$$
is the double point.  Consider the trivial family 
$$
\Cal X=\Cal K\times\Bbb C\longrightarrow U=K\times\Bbb C
$$
Assume that $X_0$ has an automorphism $\alpha=\exp \xi$ isotopic to the identity with the additional
property that its action on $H^1(X_0,\Theta_0)$ is non trivial.
 Then the family 
$$
(\exp (t\cdot \xi))_{t\in\Bbb C}
$$
defines an automorphism of $\Cal X^{\text{red}}$ which is the identity on $X_0$ and projects onto the identity of $U^{\text{red}}=\Bbb C$.
\medskip
Now the crucial point is that it also defines an automorphism $F$ of $\Cal X$ which is the identity on $X_0$. But $F$ projects onto a {\it non-trivial} morphism $f$ of $U$. Indeed, $f$ is still the identity on $U^{\text{red}}$.
But its differential is not the identity for $t\not = 0$. It may be identified with the action of $\exp (t\cdot\xi)$ on $H^1(X_t,\Theta_t)\simeq T_tU$. Hence $\Cal X$ can be obtained as a pull-back of $\Cal K$ by the map
$$
s\ :\ (0,\partial/\partial z,t)\in U=K\times\Bbb C\longmapsto (0,\partial/\partial z)\in K
$$
but also as a pull-back of $\Cal K$ by $s\circ f$. These two morphisms respect the marking at $0$ but are different, disproving the universality of $K$.  
\medskip
Observe that this argument can easily be adapted to the case where $K$ is arbitrary but non-reduced.
\medskip
Although we do not know of such a precise example (Mumford's example in \cite{Mu} has no automorphisms isotopic to the identity), it makes very plausible that Proposition 5 (especially the equivalence between (ii) and (iv)) is
not true for holomorphic families over a non-reduced base.

\medskip
Observe that Theorems 3 and 5 compared to Corollary 4 show that, surprisingly, the local isomorphism problem is fundamentally different from
the universality problem. In this last problem, there is no difference between the differentiable case and the (reduced) holomorphic case.

\vfill
\eject

\enddocument